\documentclass[11pt,a4paper]{article}

\usepackage[margin=1in]{geometry}
\usepackage{amsmath,amssymb,amsfonts,bm}
\usepackage{graphicx}
\usepackage{booktabs}
\usepackage{authblk}
\usepackage{hyperref}
\usepackage{physics}
\usepackage{mathtools}
\usepackage{xcolor}
\usepackage{xurl}
\usepackage{algorithm}
\usepackage{algorithmic}
\usepackage{placeins}
\usepackage[numbers,sort&compress]{natbib}

\hypersetup{
    colorlinks=true,
    linkcolor=blue,
    citecolor=blue,
    urlcolor=blue
}

\newcommand{\Kn}{\mathrm{Kn}}
\newcommand{\Kloc}{K}

\newcommand{\FEQ}{\mathbf{F}^{EQ}}
\newcommand{\FNS}{\mathbf{F}^{NS}}
\newcommand{\FFM}{\mathbf{F}^{FM}}
\newcommand{\FK}{\mathbf{F}^{K}}
\newcommand{\Fn}{\mathbf{F}_{n}}

\title{
A Log-Gaussian Scale-Space Limiter for Hybrid Continuum--Ballistic Gas Dynamics
}

\author[1]{Bjørn Wu}
\affil[1]{Independent Researcher, Oslo, Norway}
\date{}

\begin{document}

\maketitle

\begin{abstract}
Rarefied hypersonic gas flows involve a continuous competition between collisional relaxation and ballistic molecular transport. Classical Chapman--Enskog hydrodynamics is valid when the local Knudsen number is small, whereas free-molecular or kinetic descriptions are required when the molecular mean free path becomes comparable to or larger than the macroscopic variation length. Existing multiscale kinetic schemes provide robust treatments across these regimes but may be expensive when embedded in large-scale continuum solvers. In this work, we propose a log-Gaussian scale-space limiter for hybrid continuum--ballistic gas dynamics. The local Knudsen number is interpreted as a stochastic scale ratio, and the transition between collisional and ballistic transport is modeled as a probability partition in logarithmic Knudsen space. The resulting continuum and ballistic weights are complementary Gaussian cumulative probabilities. This construction provides super-algebraic suppression of inverse-Knudsen corrections in the continuum limit and of Chapman--Enskog corrections in the free-molecular limit. A conservative interface flux is then built by blending a Navier--Stokes--Fourier flux with a half-range Maxwellian kinetic flux. The proposed framework preserves the conservative structure of the macroscopic equations, recovers the continuum and free-molecular limits, and provides a physically interpretable transition mechanism for high-Mach rarefied flows. The method is formulated as a lightweight hybrid closure intended for future implementation in finite-volume, discrete Boltzmann, or gas-kinetic solvers. Reduced one-dimensional closure/profile comparisons against Discrete Velocity Method/Bhatnagar--Gross--Krook (DVM/BGK) Fourier and Couette data show that log-Gaussian weighting of Navier--Stokes--Fourier (NSF) and jump/slip-corrected branches improves the tested macroscopic profiles relative to NSF. The same six DVM/BGK profiles are used both as reference profiles and to calibrate \(K_0\) and \(\sigma\); therefore, the approximately \(40\%\) reduction in combined mean profile error is an in-sample calibration result for the reduced profile model rather than independent validation or a numerical validation of the proposed finite-volume face flux. Additional diagnostics assess non-equilibrium moments, internal parameter robustness, and activation of the local rarefaction indicator in analytic shock-like profiles, reduced one-velocity DVM/BGK shock layers, and two-dimensional curved shock geometries.
\end{abstract}

\noindent\textbf{Keywords:}
rarefied gas dynamics; Knudsen number; Chapman--Enskog expansion; free-molecular flow; hybrid kinetic-continuum method; log-Gaussian limiter; hypersonic flow.

\section{Introduction}

Gas flows in hypersonic and high-altitude environments frequently operate outside the strict continuum regime. In such cases, the molecular mean free path may become comparable to the macroscopic gradient length, and the local Knudsen number,
\begin{equation}
    \Kn = \frac{\lambda}{L},
\end{equation}
is no longer small. The Navier--Stokes--Fourier equations, derived from the first-order Chapman--Enskog expansion \cite{chapman_cowling,cercignani_boltzmann,cercignani_rarefied,sone,kogan,vincenti_kruger}, then lose accuracy because the distribution function is no longer a small perturbation around the local Maxwellian. In the opposite limit, when $\Kn \gg 1$, intermolecular collisions are weak over a macroscopic variation length and the flow approaches the free-molecular or ballistic transport regime.

Several families of methods have been developed to handle this multiscale problem. Direct Simulation Monte Carlo (DSMC) is reliable in the rarefied and free-molecular regimes but becomes costly in near-continuum flows \cite{bird_dsmc,bird_dsmc_book,fan_shen_ip}. Unified Gas-Kinetic Schemes (UGKS) and Discrete Unified Gas-Kinetic Schemes (DUGKS) construct interface fluxes from kinetic evolution solutions and possess asymptotic-preserving (AP) properties across a broad range of Knudsen numbers \cite{xu_huang_ugks,guo_xu_wang_dugks,guo_wang_xu_dugks_compressible,zhu_xu_ugks_decade,liu_zhu_xu_ugkwp,zhu_liu_zhong_xu_ugkwp_unstructured,xu_gks_book}. Moment methods, such as regularized 13-moment (R13) and higher-order systems, extend hydrodynamics by evolving stress and heat flux as independent non-equilibrium variables \cite{grad,struchtrup_torrilhon,torrilhon_review,levermore}. Hybrid continuum-kinetic methods couple Navier--Stokes solvers with kinetic solvers based on local rarefaction indicators \cite{boyd_breakdown,wang_boyd_breakdown,kolobov_ufs,aktas_aluru,garcia_bell_alder,wijesinghe_hadji_aluru}.

The present work proposes a complementary route. Rather than introducing a full kinetic solver everywhere, we construct a physically interpretable Knudsen-dependent probability weight that blends continuum and ballistic fluxes at the interface level. The key idea is to view the local Knudsen number not merely as a deterministic switch but as a scale-space random variable. Since $\Kn=\lambda/L$ is a ratio of positive scales, and both $\lambda$ and $L$ may vary multiplicatively in shock layers, boundary layers, and expansion regions, the logarithm of the local rarefaction measure is a natural variable. We therefore define a Gaussian distribution in $s=\ln \Kn$ and obtain continuum and ballistic weights as complementary cumulative probabilities.

The main contributions of this paper are:
\begin{enumerate}
    \raggedright
    \item A log-Gaussian scale-space interpretation of the continuum--ballistic transition.
    \item A pair of complementary Knudsen weights with correct asymptotic behavior in both $\Kn\to0$ and $\Kn\to\infty$ limits.
    \item A conservative hybrid interface flux combining Navier--Stokes--Fourier and half-range Maxwellian free-molecular fluxes.
    \item A conservative macroscopic flux formulation with a formal interpretation of selected momentum- and energy-flux corrections.
    \item A Discrete Velocity Method/Bhatnagar--Gross--Krook (DVM/BGK)-based numerical assessment for planar Fourier and Couette benchmarks,
including calibration of \(K_0\) and \(\sigma\), non-equilibrium heat-flux and shear-stress
moment diagnostics, bootstrap parameter robustness analysis, and activation tests in reduced
one-velocity DVM/BGK shock layers and two-dimensional blunt-body geometries.
\end{enumerate}

\section{Physical Motivation}

\subsection{Collision-number interpretation of the Knudsen number}

The local Knudsen number is
\begin{equation}
    \Kloc = \frac{\lambda}{L},
\end{equation}
where $\lambda$ is the molecular mean free path and $L$ is a local macroscopic variation length. Equivalently,
\begin{equation}
    \frac{1}{\Kloc} = \frac{L}{\lambda}
\end{equation}
is the expected number of mean free paths contained in one macroscopic variation length. Thus $\Kloc^{-1}$ may be interpreted as an estimate of the number of collisions available for local relaxation before a molecule traverses a region over which the macroscopic fields vary appreciably.

When
\begin{equation}
    \Kloc \ll 1,
\end{equation}
many collisions occur within the macroscopic length $L$. The gas locally relaxes toward a Maxwellian distribution, and Chapman--Enskog hydrodynamics is appropriate. When
\begin{equation}
    \Kloc \gg 1,
\end{equation}
a molecule may cross the macroscopic variation length with few or no collisions. The flow then retains upstream memory, and ballistic or kinetic transport dominates.

If the number of collisions over a distance $L$ is modeled as a Poisson process with mean
\begin{equation}
    \bar{N}=\frac{L}{\lambda}=\frac{1}{\Kloc},
\end{equation}
then the probability of no collision over the macroscopic scale is
\begin{equation}
    P_0 = \exp\left(-\frac{1}{\Kloc}\right),
\end{equation}
and the probability of at least one collision is
\begin{equation}
    P_{\ge 1} = 1-\exp\left(-\frac{1}{\Kloc}\right).
\end{equation}
This elementary argument suggests that continuum and ballistic transport may be interpreted probabilistically rather than as a sharp deterministic switch. It should not, however, be read as a derivation of the final log-Gaussian limiter. The Poisson no-collision probability and the log-Gaussian ballistic weight are distinct stochastic transition models. If
\begin{equation}
    s=\ln\left(\frac{K}{K_0}\right),
\end{equation}
then $K=K_0 e^s$, and the Poisson expression can be rewritten as
\begin{equation}
    P_0(K)=\exp\left(-\frac{1}{K}\right)
    =\exp\left[-\frac{\exp(-s)}{K_0}\right].
\end{equation}
Thus, as a function of $s$, the Poisson no-collision model has the form of a shifted unit-scale Gumbel cumulative distribution in log-Knudsen space, with location parameter $-\ln K_0$, whereas the proposed limiter uses a Gaussian cumulative distribution in the same log-Knudsen variable. Both are smooth maps from $K>0$ to $(0,1)$, but they have different shapes, tail behavior, and modeling roles. In the present paper, the Poisson expression is used only to motivate probabilistic collision/ballistic language. The final log-Gaussian form is instead chosen to provide a symmetric tunable transition in log-scale space and super-algebraic endpoint suppression of invalid asymptotic branches.

\subsection{Scale-space uncertainty and log-Gaussian transition}

In hypersonic rarefied flows, neither $\lambda$ nor $L$ is a fixed deterministic quantity. The mean free path depends on local density, temperature, and molecular interaction model, while the gradient length depends on which field is used:
\begin{equation}
    L_\rho=\frac{\rho}{|\nabla\rho|},
    \qquad
    L_T=\frac{T}{|\nabla T|},
    \qquad
    L_u=\frac{|\mathbf{u}|+a}{\|\nabla\mathbf{u}\|}.
\end{equation}
Moreover, shock layers, wall Knudsen layers, and expansion fans introduce directional and history-dependent non-equilibrium effects. Therefore, the critical transition between collisional and ballistic behavior should not be treated as a single deterministic threshold.

Because $\Kloc=\lambda/L$ is a ratio of positive scales, a multiplicative uncertainty model suggests that $\ln \Kloc$ is the natural transition variable. We introduce
\begin{equation}
    s = \ln\left(\frac{\Kloc}{K_0}\right),
\end{equation}
where $K_0$ is the nominal transition Knudsen number. The parameter $K_0^{-1}$ can be interpreted as the critical number of collisions over a macroscopic scale. If $K_0=1$, the transition is centered at one collision per macroscopic variation length. If $K_0=0.1$, the ballistic correction is activated when fewer than approximately ten collisions occur across $L$.

We assume that the critical transition scale in log-Knudsen space is Gaussian distributed with standard deviation $\sigma$. The continuum and ballistic weights are then defined as complementary cumulative probabilities:
\begin{equation}
    W_c(\Kloc)
    =
    \frac{1}{2}
    \operatorname{erfc}\left[
    \frac{\ln(\Kloc/K_0)}{\sqrt{2}\sigma}
    \right],
    \label{eq:Wc}
\end{equation}
\begin{equation}
    W_f(\Kloc)
    =
    1-W_c(\Kloc)
    =
    \frac{1}{2}
    \operatorname{erfc}\left[
    -\frac{\ln(\Kloc/K_0)}{\sqrt{2}\sigma}
    \right].
    \label{eq:Wf}
\end{equation}
Here $W_c$ is interpreted as a modeled probability weight for local collisional dominance, and $W_f$ is interpreted as a modeled probability weight for local ballistic dominance.

\section{Asymptotic Properties of the Log-Gaussian Weight}

The proposed weights satisfy
\begin{equation}
    \Kloc\to0:
    \qquad
    W_c\to1,
    \qquad
    W_f\to0,
\end{equation}
and
\begin{equation}
    \Kloc\to\infty:
    \qquad
    W_c\to0,
    \qquad
    W_f\to1.
\end{equation}
More importantly, the log-Gaussian tails suppress divergent asymptotic branches. As $\Kloc\to0$, let
\begin{equation}
    y=-\ln \Kloc \to \infty.
\end{equation}
Then $W_f$ has the standard Gaussian-tail asymptotic form
\begin{equation}
    W_f(\Kloc)
    \sim
    \frac{\sigma}{\sqrt{2\pi}\,|\ln(\Kloc/K_0)|}
    \exp\left[
    -\frac{\ln^2(\Kloc/K_0)}{2\sigma^2}
    \right].
\end{equation}
For any finite $m>0$,
\begin{equation}
    \log\left(W_f(\Kloc)\Kloc^{-m}\right)
    =
    -\frac{(y+\ln K_0)^2}{2\sigma^2}
    +my
    +O(\ln y)
    \to -\infty.
\end{equation}
Thus inverse-Knudsen terms are super-algebraically suppressed in the continuum limit. Similarly, as $\Kloc\to\infty$,
\begin{equation}
    W_c(\Kloc)\Kloc^n\to0
\end{equation}
for any finite $n>0$. Therefore, the weight can consistently combine positive powers and inverse powers of $\Kloc$ without introducing endpoint divergences.

This motivates the formal double-asymptotic representation
\begin{equation}
    f
    =
    f^M
    +
    W_c(\Kloc)\sum_{n=1}^{N}\Kloc^n f_c^{(n)}
    +
    W_f(\Kloc)\sum_{m=1}^{M}\Kloc^{-m}f_f^{(m)},
    \label{eq:double_expansion}
\end{equation}
where $f^M$ is the local Maxwellian, $f_c^{(n)}$ are continuum-side Chapman--Enskog corrections, and $f_f^{(m)}$ are ballistic-side inverse-Knudsen corrections. The log-Gaussian weights ensure that each branch is active only in the regime where it is asymptotically meaningful.

\section{Local Rarefaction Indicator}

The local rarefaction parameter is defined as
\begin{equation}
    \Kloc = \max(K_\rho,K_T,K_u),
    \label{eq:Klocal}
\end{equation}
where
\begin{equation}
    K_\rho=\lambda\frac{|\nabla\rho|}{\rho},
\end{equation}
\begin{equation}
    K_T=\lambda\frac{|\nabla T|}{T},
\end{equation}
\begin{equation}
    K_u=\lambda\frac{\|\nabla\mathbf{u}\|}{|\mathbf{u}|+a},
\end{equation}
and
\begin{equation}
    a=\sqrt{\gamma RT}
\end{equation}
is the local sound speed. Alternative definitions may include pressure gradients, heat-flux-based indicators, or thermodynamic non-equilibrium measures. In a discrete Boltzmann implementation, one may define
\begin{equation}
    K_{\mathrm{TNE}}
    =
    C_2\frac{\|\Delta_2^*\|}{p}
    +
    C_3\frac{\|\Delta_{3,1}^*\|}{p\sqrt{RT}}
    +\cdots,
\end{equation}
where $\Delta_2^*$ and $\Delta_{3,1}^*$ are non-equilibrium stress and heat-flux moments. Such a definition links the weight directly to the departure of the local distribution function from the hydrodynamic slow manifold.

For a finite-volume face, we use
\begin{equation}
    \Kloc_{face}=\max(\Kloc_L,\Kloc_R),
    \label{eq:Kface}
\end{equation}
which biases the flux toward kinetic treatment whenever either side of the interface is locally rarefied.

\section{Conservative Hybrid Flux}

We consider the conservative variables
\begin{equation}
    \mathbf{U}
    =
    \begin{bmatrix}
    \rho\\
    \rho u_x\\
    \rho u_y\\
    E
    \end{bmatrix},
\end{equation}
with
\begin{equation}
    E=\rho c_vT+\frac{1}{2}\rho|\mathbf{u}|^2,
\end{equation}
and the ideal-gas equation of state
\begin{equation}
    p=\rho RT.
\end{equation}

At a cell interface with unit normal $\mathbf{n}$, the proposed hybrid flux is
\begin{equation}
    \Fn
    =
    W_c(\Kloc_{face})\FNS_n
    +
    W_f(\Kloc_{face})\FK_n.
    \label{eq:hybrid_flux}
\end{equation}
Here $\FNS_n$ is the Navier--Stokes--Fourier flux and $\FK_n$ is a kinetic or ballistic flux. The interface flux is shared by neighboring cells with opposite sign, so conservation follows immediately in a finite-volume discretization.

This conservation statement is deliberately limited. A single shared face flux gives discrete conservation of mass, momentum, and energy, but it does not by itself imply positivity preservation of density or temperature, entropy stability, monotonicity, or nonlinear stability. These stronger properties require additional limiters, timestep restrictions, or separate analysis and are not claimed in the present work. The present contribution establishes a conservative and limiting-consistent flux construction; positivity-preserving, entropy-stable, and AP analyses are left for future work.

\section{Continuum Flux}

The Navier--Stokes--Fourier flux in the normal direction is
\begin{equation}
    \FNS_n =
    \begin{bmatrix}
    \rho u_n\\
    \rho u_n\mathbf{u}+p\mathbf{n}-\boldsymbol{\sigma}^{NS}\cdot\mathbf{n}\\
    (E+p)u_n-(\boldsymbol{\sigma}^{NS}\cdot\mathbf{u})\cdot\mathbf{n}+q_n^{NS}
    \end{bmatrix},
    \label{eq:NS_flux}
\end{equation}
where
\begin{equation}
    u_n=\mathbf{u}\cdot\mathbf{n}.
\end{equation}
The viscous stress is
\begin{equation}
    \sigma_{ij}^{NS}
    =
    \mu
    \left[
    \partial_i u_j+\partial_j u_i
    -
    \frac{2}{D}(\nabla\cdot\mathbf{u})\delta_{ij}
    \right]
    +
    \zeta(\nabla\cdot\mathbf{u})\delta_{ij},
    \label{eq:NS_stress}
\end{equation}
and the Fourier heat flux is
\begin{equation}
    \mathbf{q}^{NS}=-\kappa\nabla T.
    \label{eq:Fourier}
\end{equation}
The heat conductivity is
\begin{equation}
    \kappa=\frac{\mu c_p}{Pr}.
\end{equation}
For a monatomic gas, $\zeta$ may be set to zero. For hypersonic air or gases with internal degrees of freedom, bulk viscosity and temperature-dependent transport coefficients may be included.

\section{Half-Range Maxwellian Free-Molecular Flux}

The ballistic component is constructed from a half-range kinetic flux. Let $\boldsymbol{\xi}$ be the molecular velocity and
\begin{equation}
    \xi_n=\boldsymbol{\xi}\cdot\mathbf{n}.
\end{equation}
Define the collision invariants
\begin{equation}
    \boldsymbol{\psi}(\boldsymbol{\xi},\eta)
    =
    \begin{bmatrix}
    1\\
    \xi_x\\
    \xi_y\\
    \frac{1}{2}(\xi^2+\eta^2)
    \end{bmatrix},
    \label{eq:psi}
\end{equation}
where $\eta$ denotes internal degrees of freedom. The free-molecular flux is
\begin{equation}
    \FFM_n
    =
    \int_{\xi_n>0}
    \xi_n\boldsymbol{\psi}f_L^M\,\dd\boldsymbol{\xi}\dd\eta
    +
    \int_{\xi_n<0}
    \xi_n\boldsymbol{\psi}f_R^M\,\dd\boldsymbol{\xi}\dd\eta.
    \label{eq:FM_flux_integral}
\end{equation}
The Maxwellian with internal degrees of freedom is
\begin{equation}
    f^M
    =
    \rho
    \frac{1}{(2\pi RT)^{D/2}}
    \frac{1}{(2\pi RT)^{b/2}}
    \exp\left[
    -
    \frac{
    |\boldsymbol{\xi}-\mathbf{u}|^2+\eta^2
    }{2RT}
    \right],
    \label{eq:Maxwellian}
\end{equation}
where
\begin{equation}
    \gamma=\frac{D+b+2}{D+b},
\end{equation}
or equivalently,
\begin{equation}
    b=\frac{2}{\gamma-1}-D.
\end{equation}

For a face-local coordinate system with normal velocity $u_n$, tangential velocity $u_t$, and
\begin{equation}
    \theta=RT,
\end{equation}
define
\begin{equation}
    a=\frac{u_n}{\sqrt{2\theta}},
\end{equation}
\begin{equation}
    A^+=\frac{1}{2}\operatorname{erfc}(-a),
    \qquad
    A^-=\frac{1}{2}\operatorname{erfc}(a),
\end{equation}
and
\begin{equation}
    B=\sqrt{\frac{\theta}{2\pi}}\exp(-a^2).
\end{equation}

The positive half-range Maxwellian moments are
\begin{align}
    M_0^+ &= \rho A^+,\\
    M_1^+ &= \rho(u_nA^+ + B),\\
    M_2^+ &= \rho[(u_n^2+\theta)A^+ + u_nB],\\
    M_3^+ &= \rho[(u_n^3+3u_n\theta)A^+ +(u_n^2+2\theta)B].
\end{align}
The negative half-range moments are
\begin{align}
    M_0^- &= \rho A^-,\\
    M_1^- &= \rho(u_nA^- - B),\\
    M_2^- &= \rho[(u_n^2+\theta)A^- - u_nB],\\
    M_3^- &= \rho[(u_n^3+3u_n\theta)A^- -(u_n^2+2\theta)B].
\end{align}

For a two-dimensional gas with one tangential velocity, the half-range flux contributions are
\begin{align}
    F_\rho^\pm &= M_1^\pm,\\
    F_{\rho u_n}^\pm &= M_2^\pm,\\
    F_{\rho u_t}^\pm &= u_tM_1^\pm,\\
    F_E^\pm
    &=
    \frac{1}{2}M_3^\pm
    +
    \frac{1}{2}
    \left[
    u_t^2+(1+b)\theta
    \right]M_1^\pm.
\end{align}
The total free-molecular interface flux is
\begin{equation}
    \FFM_n
    =
    \mathbf{F}^+(\rho_L,u_{n,L},u_{t,L},T_L)
    +
    \mathbf{F}^-(\rho_R,u_{n,R},u_{t,R},T_R).
    \label{eq:FM_sum}
\end{equation}
The local momentum flux is then rotated back from $(n,t)$ coordinates to Cartesian coordinates.

\section{BGK-Corrected Kinetic Flux}

The pure free-molecular flux is appropriate in the limit $\Kloc\gg1$, but in the transition regime a finite collision correction is desirable. A simple BGK-inspired kinetic flux may be written as \cite{bgk,shakhov,holway}
\begin{equation}
    \FK_n
    =
    \alpha\FFM_n
    +
    (1-\alpha)\FEQ_n,
    \label{eq:BGK_flux}
\end{equation}
where $\FEQ_n$ is the equilibrium Euler flux and
\begin{equation}
    \alpha
    =
    \frac{\tau}{\Delta t}
    \left(
    1-e^{-\Delta t/\tau}
    \right).
    \label{eq:alpha}
\end{equation}
When $\tau\gg\Delta t$,
\begin{equation}
    \alpha\to1,
\end{equation}
and the kinetic flux approaches the free-molecular flux. When $\tau\ll\Delta t$,
\begin{equation}
    \alpha\to0,
\end{equation}
and the kinetic flux approaches the equilibrium flux.

The final proposed interface flux is therefore
\begin{equation}
    \boxed{
    \Fn
    =
    W_c(\Kloc_{face})\FNS_n
    +
    W_f(\Kloc_{face})
    \left[
    \alpha\FFM_n
    +
    (1-\alpha)\FEQ_n
    \right].
    }
    \label{eq:final_flux}
\end{equation}
This expression is conservative, inexpensive, and naturally recovers the continuum and ballistic limits.

\section{Macroscopic Equations}

The finite-volume method is defined directly by the shared numerical face flux. In formal continuum notation this corresponds to the macroscopic conservation law
\begin{equation}
    \partial_t\mathbf{U}
    +
    \nabla\cdot\mathbf{F}^{hyb}
    =
    0,
    \label{eq:macro_conservation}
\end{equation}
where the Knudsen-weighted hybrid flux is the finite-volume flux
\begin{equation}
    \mathbf{F}^{hyb}
    =
    W_c\mathbf{F}^{NS}
    +
    W_f\mathbf{F}^{K}.
    \label{eq:macro_hybrid_flux}
\end{equation}
This conservative flux formulation is the primary macroscopic statement of the method. The kinetic contribution generally modifies all conservative flux components, including the mass flux. It therefore need not admit an exact decomposition into the standard convective flux plus only effective stress and heat-flux terms. Selected momentum- and energy-flux corrections may be interpreted formally in terms of effective stress and heat-flux contributions, but this interpretation is not an equivalent local PDE representation of the complete hybrid flux. In the rarefied limit the model should be understood as a conservative macroscopic balance equipped with a nonlocal kinetic flux closure.

\section{Limiting Behavior}

\subsection{Continuum limit}

When
\begin{equation}
    \Kloc\to0,
\end{equation}
we have
\begin{equation}
    W_c\to1,
    \qquad
    W_f\to0.
\end{equation}
Thus
\begin{equation}
    \Fn\to\FNS_n,
\end{equation}
and the model recovers the Navier--Stokes--Fourier equations.

\subsection{Free-molecular limit}

When
\begin{equation}
    \Kloc\to\infty,
\end{equation}
we have
\begin{equation}
    W_c\to0,
    \qquad
    W_f\to1.
\end{equation}
If $\tau\gg\Delta t$, then
\begin{equation}
    \alpha\to1,
\end{equation}
and
\begin{equation}
    \Fn\to\FFM_n.
\end{equation}
Therefore, the model recovers the half-range Maxwellian free-molecular flux.

\subsection{Transition regime}

When
\begin{equation}
    \Kloc=O(1),
\end{equation}
both weights are active. The flux is a probability-weighted combination of continuum relaxation and ballistic transport. The parameter $\sigma$ controls the width of the transition region in log-Knudsen space.

\section{Numerical Implementation}

A prospective finite-volume implementation would proceed as follows.

\begin{algorithm}[htbp]
\caption{Log-Gaussian hybrid continuum--ballistic flux}
\begin{algorithmic}[1]
\STATE Reconstruct left and right interface states:
\[
(\rho_L,\mathbf{u}_L,T_L),\qquad
(\rho_R,\mathbf{u}_R,T_R).
\]
\STATE Compute local mean free path $\lambda$ and local rarefaction indicators $K_\rho,K_T,K_u$.
\STATE Set
\[
K_{face}=\max(K_L,K_R).
\]
\STATE Evaluate
\[
W_c(K_{face}),\qquad W_f(K_{face}).
\]
\STATE Compute the Navier--Stokes--Fourier flux $\FNS_n$.
\STATE Compute the half-range Maxwellian free-molecular flux $\FFM_n$ using analytic error-function moments.
\STATE Compute the BGK relaxation coefficient
\[
\alpha=
\frac{\tau}{\Delta t}
\left(1-e^{-\Delta t/\tau}\right).
\]
\STATE Construct
\[
\Fn
=
W_c\FNS_n
+
W_f
\left[
\alpha\FFM_n
+
(1-\alpha)\FEQ_n
\right].
\]
\STATE Update cell averages conservatively.
\end{algorithmic}
\end{algorithm}

The proposed formulation would require no global domain decomposition between continuum and kinetic regions. The transition would remain local, smooth, and controlled by physically interpretable scale-space probabilities in a future solver implementation.


\FloatBarrier

\section{Numerical Experiments and Diagnostics}

This section tests the proposed log-Gaussian continuum--ballistic limiter and calibrates its transition parameters against one-dimensional Discrete Velocity Method/Bhatnagar--Gross--Krook (DVM/BGK) reference solutions. The experiments focus on five questions: whether the weight behaves smoothly across
Knudsen regimes, whether rarefaction corrections improve planar wall-bounded profiles
relative to NSF, how the parameters \(K_0\) and \(\sigma\) should be chosen, whether the
local indicator activates in shock-like high-gradient layers, and whether the calibrated
parameters show basic holdout robustness. In the planar benchmark problems below, the global Knudsen number is used as the controlled rarefaction parameter; in multidimensional finite-volume applications it is replaced by the local face indicator \(K_{\mathrm{face}}\).

\subsection{Numerical protocol and status of the comparison models}
\label{sec:numerical_protocol}

The planar Fourier and Couette tests in Sections~12.3--12.5 are closure/profile diagnostics against deterministic DVM/BGK reference data, not demonstrations of a fully time-dependent multidimensional production finite-volume solver. The same six DVM/BGK profiles are used both as reference profiles for the reduced comparisons and as the calibration dataset for \(K_0\) and \(\sigma\). They therefore support development and internal assessment of the reduced profile model, but they do not constitute independent validation across different flow classes or external datasets. The DVM/BGK boundary-value problem is solved as the kinetic reference. The NSF, hard-switch, algebraic, jump/slip, and log-Gaussian curves are reduced one-dimensional closure predictions evaluated on the same geometry, wall data, Knudsen numbers, and error metric. The conservative flux formula in Sections~5--8 is the proposed finite-volume interface form; these reduced profile comparisons do not numerically validate that face flux in a complete finite-volume solver. Full multidimensional finite-volume solver implementation, time advancement, and validation are reserved for future work.

For the DVM/BGK references, the wall-normal coordinate is discretized on a uniform one-dimensional grid and iterated to a steady state with diffuse wall boundary conditions. The walls prescribe temperature in the Fourier problem and wall velocity plus temperature in the Couette problem. The mean free path is set by the prescribed global Knudsen number for each planar benchmark case, and the same global Knudsen number is used as the controlled rarefaction parameter for the reduced comparison curves. In a multidimensional implementation, this controlled global value is replaced by the local face value \(K_{\mathrm{face}}\) defined in Eq.~\eqref{eq:Kface}.

The DVM/BGK reference scripts use a uniform physical grid with \(N_y=121\) wall-normal points and a tensor-product velocity grid with \(33\times33\) velocity nodes on \([-7,7]^2\). The relaxation time is set nondimensionally as \(\tau=\max(Kn,10^{-6})\). A first-order upwind sweep is used for positive and negative wall-normal molecular velocities, near-zero wall-normal velocities are locally relaxed, and the iteration stops when the relative distribution-function change falls below \(3.0\times10^{-8}\), with a maximum of 2500 iterations. These settings define the reference DVM/BGK boundary-value solves used for the profile and moment diagnostics. They are separate from the proposed multidimensional finite-volume face-flux formulation, for which production-solver reconstruction, time integration, and solver-level convergence studies are left for future work.

The comparison branches are defined as follows. The NSF branch is the continuum no-slip/no-jump prediction obtained from the Fourier or Couette continuum solution. The hard-switch branch uses
\begin{equation}
W_f^{hard}(K)=
\begin{cases}
0, & K<K_0,\\
1, & K\ge K_0,
\end{cases}
\qquad W_c^{hard}=1-W_f^{hard}.
\end{equation}
The algebraic branch uses
\begin{equation}
W_f^{alg}(K)=\frac{K^p}{K^p+K_0^p},\qquad W_c^{alg}=1-W_f^{alg},
\end{equation}
with the same transition center used for the corresponding comparison. In the branch-comparison figures, the hard-switch and algebraic comparisons use \(K_0=0.1\), and the algebraic exponent is \(p=2\), matching the comparison scripts. The log-Gaussian branch uses the continuum and ballistic weights in Eqs.~\eqref{eq:Wc}--\eqref{eq:Wf}. The jump branch in the Fourier benchmark applies a rarefaction-corrected temperature-jump profile relative to the DVM/BGK wall data, while the slip branch in the Couette benchmark applies a rarefaction-corrected velocity-slip profile. These branches are included as reduced closure comparators rather than as independent kinetic solvers.

For the reduced comparison branches, the predicted macroscopic profile is constructed at the profile level as
\begin{equation}
\phi_{model}(y;K)=W_c(K)\phi_{NSF}(y)+W_f(K)\phi_R(y),
\label{eq:reduced_profile_blend}
\end{equation}
where \(\phi=T\) for Fourier flow and \(\phi=u\) for Couette flow. The rarefied profile \(\phi_R\) denotes the temperature-jump-corrected Fourier profile in the heat-transfer benchmark and the velocity-slip-corrected Couette profile in the shear-flow benchmark. Specifically, with the normalized wall-normal coordinate \(y\in[0,1]\), the Fourier comparison uses
\begin{align}
T_{NSF}(y)&=T_c+(T_h-T_c)y,\\
T_R(y;K)&=T_c+(T_h-T_c)\frac{y+B_TK}{1+2B_TK},
\end{align}
with \(T_c=1\), \(T_h=1.5\), and \(B_T=2.18\). The Couette comparison uses
\begin{align}
u_{NSF}(y)&=U_wy,\\
u_R(y;K)&=U_w\frac{y+A_uK}{1+2A_uK},
\end{align}
with \(U_w=0.2\) and \(A_u=1.146\). The constants \(B_T=2.18\) and \(A_u=1.146\) are prescribed comparator coefficients in the reduced jump/slip branches. They are held fixed throughout the log-Gaussian calibration and are not inferred from the six DVM/BGK profiles used to calibrate \(K_0\) and \(\sigma\). They are not proposed as new universal gas-surface coefficients in this work; their role is only to define reproducible reduced comparison branches, not to modify the proposed log-Gaussian weighting rule. The hard-switch, algebraic, and log-Gaussian comparison curves differ through the choice of \(W_f(K)\) and \(W_c(K)=1-W_f(K)\) in Eq.~\eqref{eq:reduced_profile_blend}. The pure jump and slip curves correspond to using \(\phi_R\) directly. This construction is a reduced profile diagnostic; it is not an additional kinetic boundary-value solve.

The profile error is computed as
\begin{equation}
E_{\phi}=\frac{\|\phi_{model}-\phi_{DVM}\|_2}{\|\phi_{DVM}\|_2},
\end{equation}
where \(\phi=T\) for Fourier flow and \(\phi=u\) for Couette flow. The same discrete wall-normal grid is used for the model and reference profiles when computing this metric. This distinction specifies which equations are solved as kinetic references and which curves are reduced closure predictions.

\subsection{Weight behavior}

The ballistic weight is
\[
W_f(K)=
\frac12\operatorname{erfc}
\left(
-\frac{\ln(K/K_0)}{\sqrt{2}\sigma}
\right),
\qquad
W_c(K)=1-W_f(K).
\]
Here \(K_0\) sets the transition center in logarithmic Knudsen space and \(\sigma\) sets the transition width. The default setting is \(K_0=0.1,\sigma=1.0\).

\begin{figure}[htbp]
\centering
\includegraphics[width=0.70\textwidth]{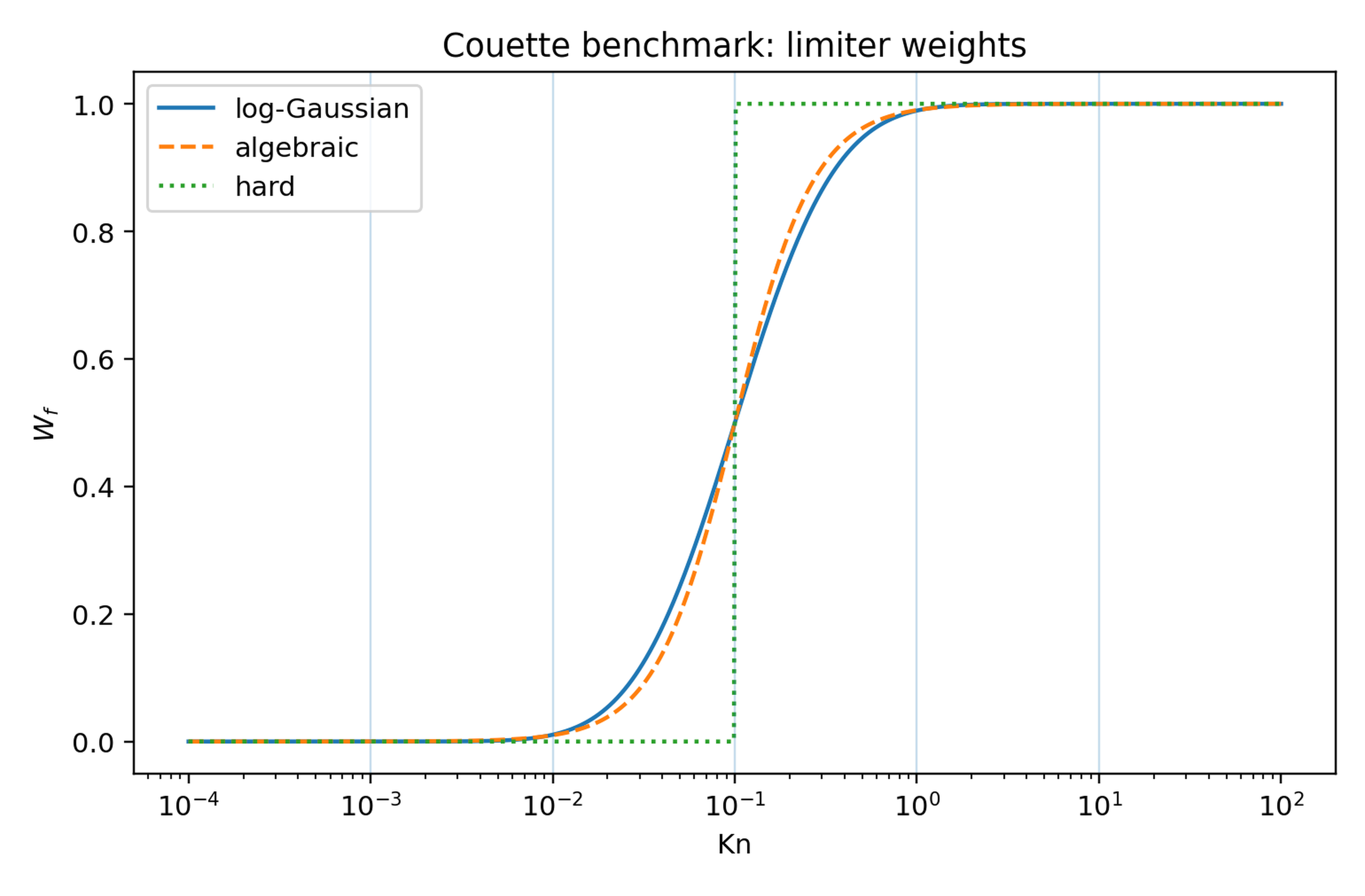}
\caption{Comparison of log-Gaussian, algebraic, and hard-switch ballistic weights.}
\label{fig:weight_curves}
\end{figure}

\subsection{DVM/BGK Fourier benchmark}

\begin{figure}[htbp]
\centering
\includegraphics[width=\textwidth]{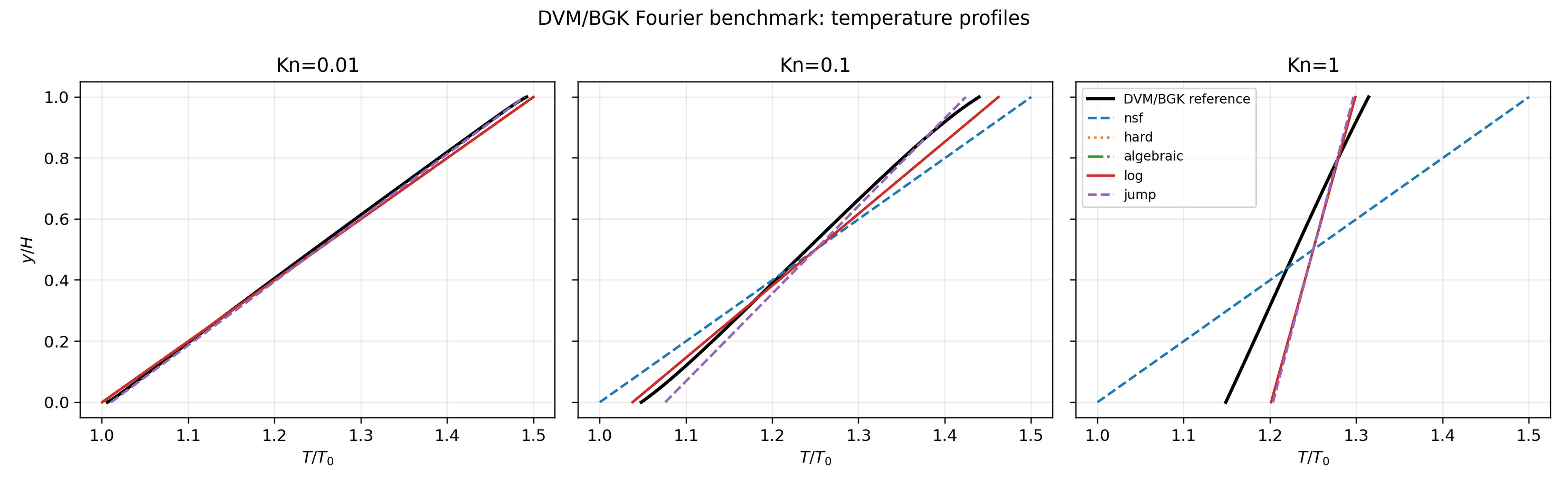}
\caption{DVM/BGK Fourier benchmark. Temperature profiles from the kinetic reference are compared with continuum and rarefied-branch predictions.}
\label{fig:dvm_fourier_profiles}
\end{figure}

Fourier heat transfer is considered between two diffuse walls \cite{cercignani_lampis,sharipov_seleznev} with \(T_c=1\) and \(T_h=1.5\). A stationary DVM/BGK equation,
\[
v_y\partial_y f=\frac{M[f]-f}{\tau},
\]
is solved in the wall-normal coordinate and used as a kinetic reference. As clarified in Section~\ref{sec:numerical_protocol}, Figure~\ref{fig:dvm_fourier_profiles} compares the DVM/BGK temperature profiles with reduced NSF, hard-switch, algebraic, log-Gaussian, and jump-branch closure predictions.

\subsection{DVM/BGK Couette benchmark}

\begin{figure}[htbp]
\centering
\includegraphics[width=\textwidth]{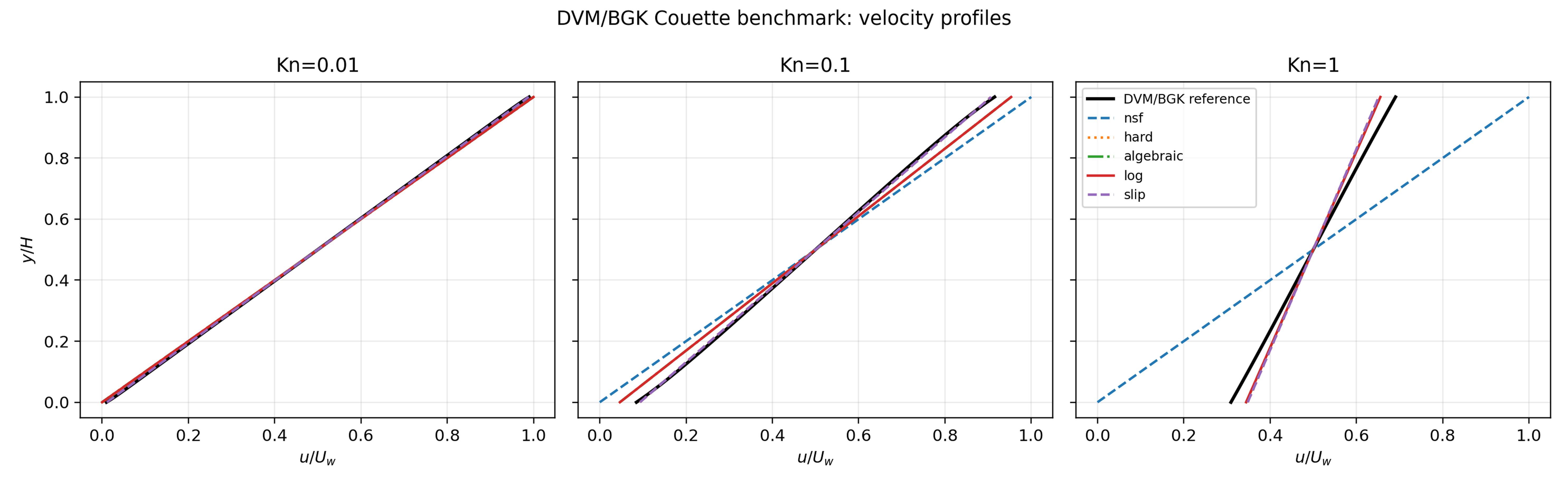}
\caption{DVM/BGK Couette benchmark. Velocity profiles from the kinetic reference are compared with continuum and rarefied-branch predictions.}
\label{fig:dvm_couette_profiles}
\end{figure}

Couette flow is considered between diffuse walls \cite{cercignani_lampis,sharipov_seleznev}. The lower wall is stationary and the upper wall moves with \(U_w=0.2\), with fixed wall temperature \(T_w=1\). Figure~\ref{fig:dvm_couette_profiles} shows that the NSF no-slip profile increasingly deviates from the DVM/BGK reference as \(Kn\) grows, while reduced slip-type and log-Gaussian closure predictions capture the rarefaction-induced velocity-slip trend.

\subsection{Parameter calibration}

The default pair \(K_0=0.1,\sigma=1.0\) is physically transparent but not necessarily optimal for wall-bounded rarefied flows. We therefore scan \(K_0\) and \(\sigma\) and compute
\[
E_\phi=
\frac{\|\phi_{\mathrm{model}}-\phi_{\mathrm{DVM}}\|_2}
{\|\phi_{\mathrm{DVM}}\|_2},
\]
where \(\phi=T\) for Fourier flow and \(\phi=u\) for Couette flow. Within the scanned range, the best tested pair is
\[
K_0=0.03,\qquad \sigma=2.5.
\]

\begin{figure}[htbp]
\centering
\includegraphics[width=0.78\textwidth]{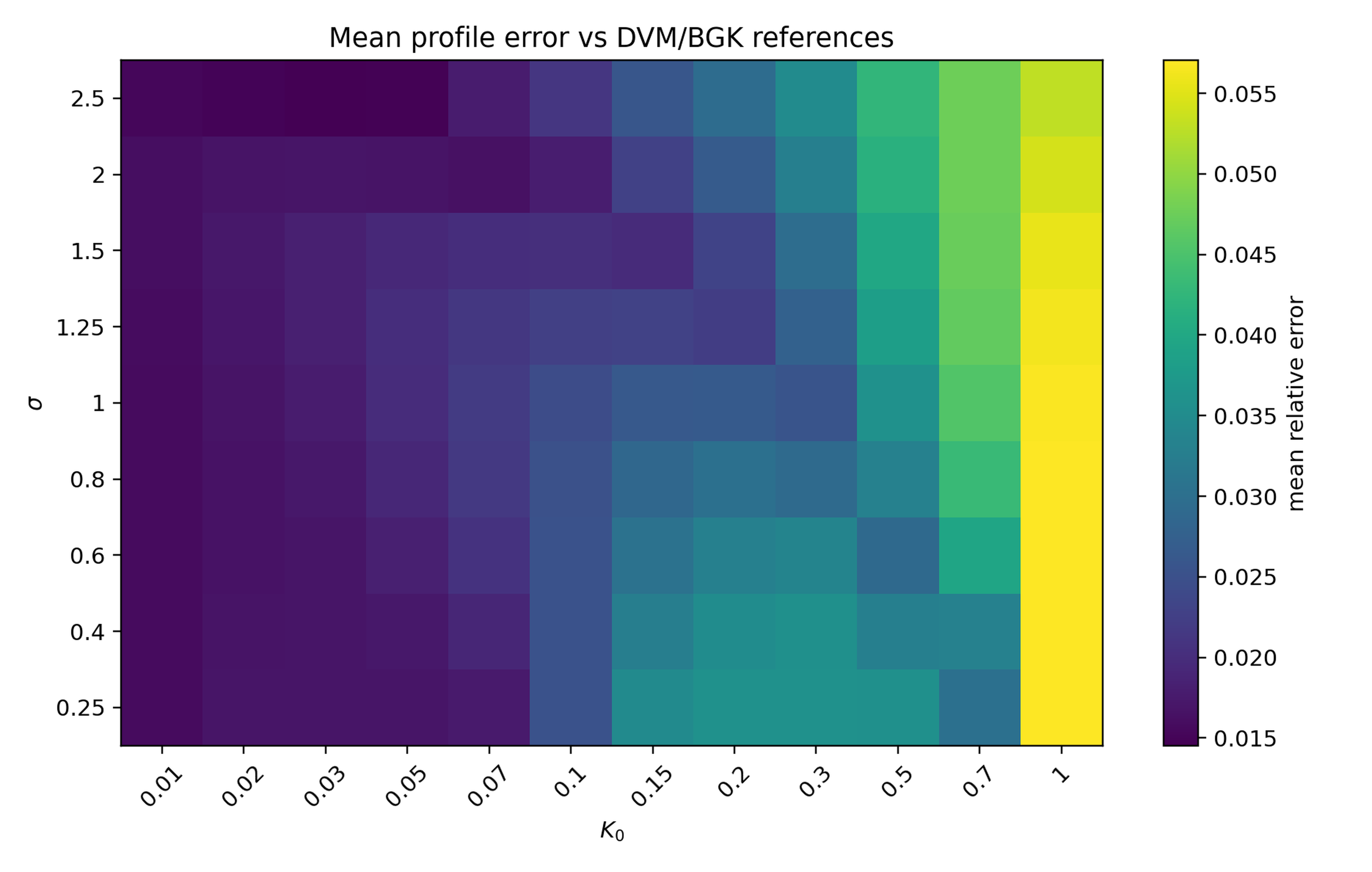}
\caption{Mean profile error against DVM/BGK Fourier and Couette references as a function of \(K_0\) and \(\sigma\).}
\label{fig:parameter_fit_heatmap}
\end{figure}

Table~\ref{tab:default_vs_fitted} summarizes the default and calibrated choices. The calibrated pair reduces the combined mean profile error from \(2.44\times10^{-2}\) to \(1.45\times10^{-2}\), approximately a \(40\%\) reduction, on the same six Fourier/Couette DVM/BGK cases used to select \(K_0\) and \(\sigma\). This reduction should therefore be interpreted as an in-sample calibration result for the reduced profile model, not as independent validation.

\begin{table}[htbp]
\centering
\caption{Default and DVM/BGK-calibrated log-Gaussian parameters. The calibrated result is in-sample with respect to the six Fourier/Couette DVM/BGK profiles used for parameter selection.}
\label{tab:default_vs_fitted}
\begin{tabular}{llrrrr}
\hline
Choice & Parameters & Fourier mean & Couette mean & Overall mean & Max error \\
\hline
default & $K_0=0.1,\quad \sigma=1.0$ & 1.3819e-02 & 3.4916e-02 & 2.4368e-02 & 5.3505e-02 \\
fitted & $K_0=0.03,\quad \sigma=2.5$ & 1.0526e-02 & 1.8540e-02 & 1.4533e-02 & 3.5527e-02 \\
\hline
\end{tabular}
\end{table}

\begin{figure}[htbp]
\centering
\includegraphics[width=\textwidth]{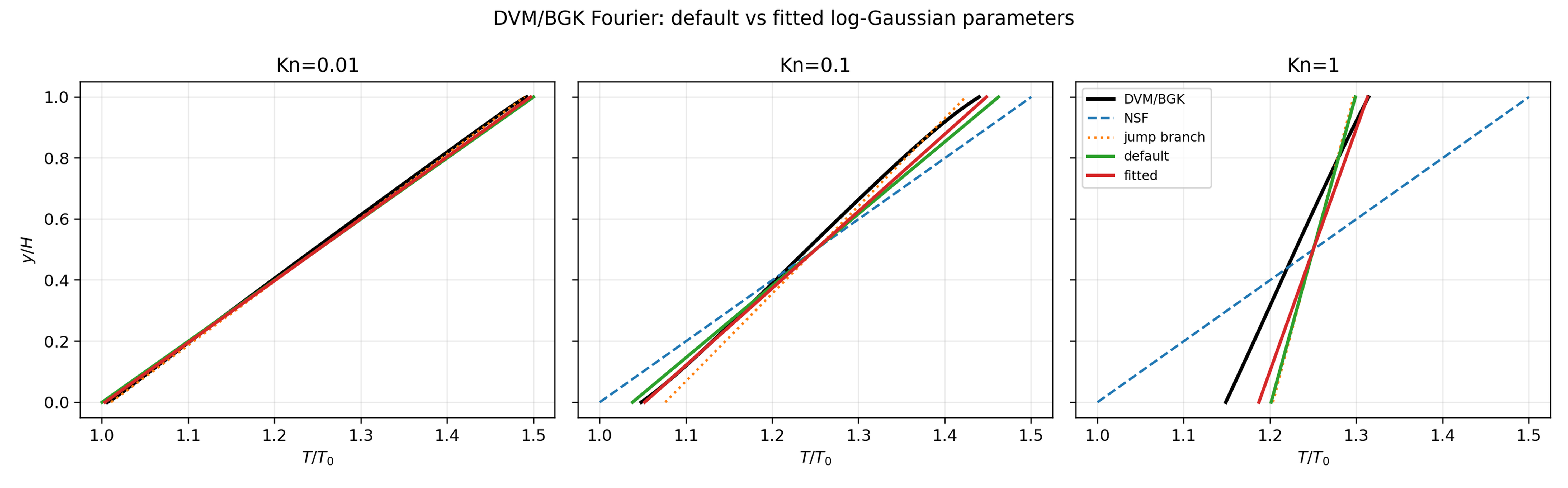}
\caption{DVM/BGK Fourier benchmark: comparison of default and calibrated log-Gaussian parameters.}
\label{fig:fourier_default_fitted}
\end{figure}

\begin{figure}[htbp]
\centering
\includegraphics[width=\textwidth]{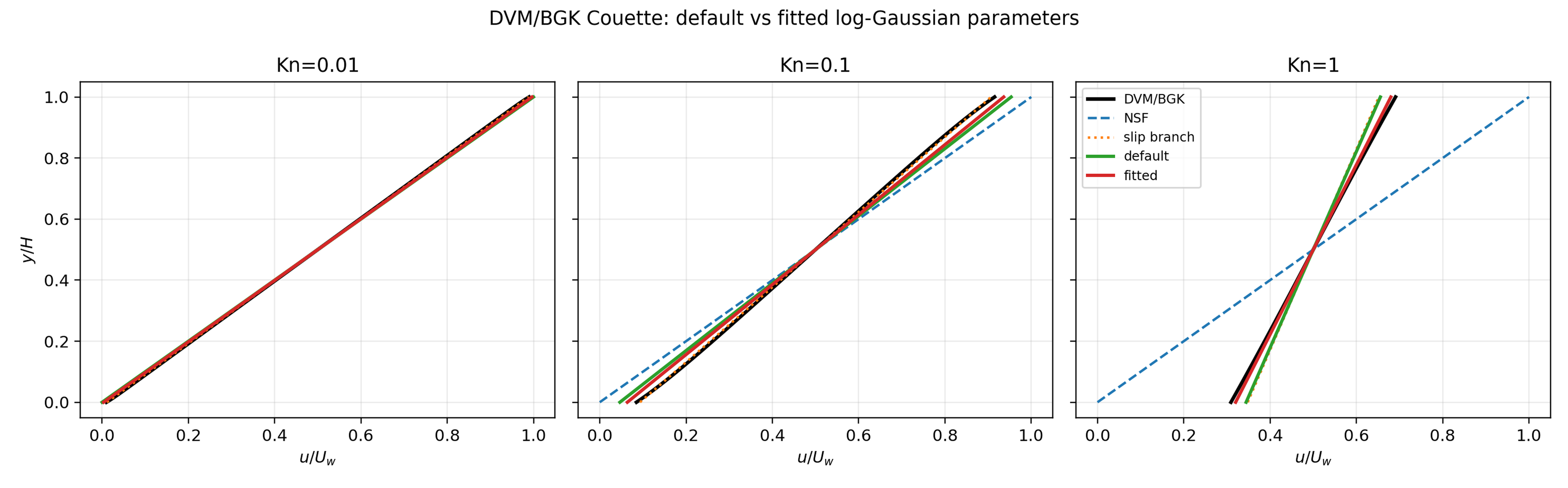}
\caption{DVM/BGK Couette benchmark: comparison of default and calibrated log-Gaussian parameters.}
\label{fig:couette_default_fitted}
\end{figure}

\begin{figure}[htbp]
\centering
\includegraphics[width=\textwidth]{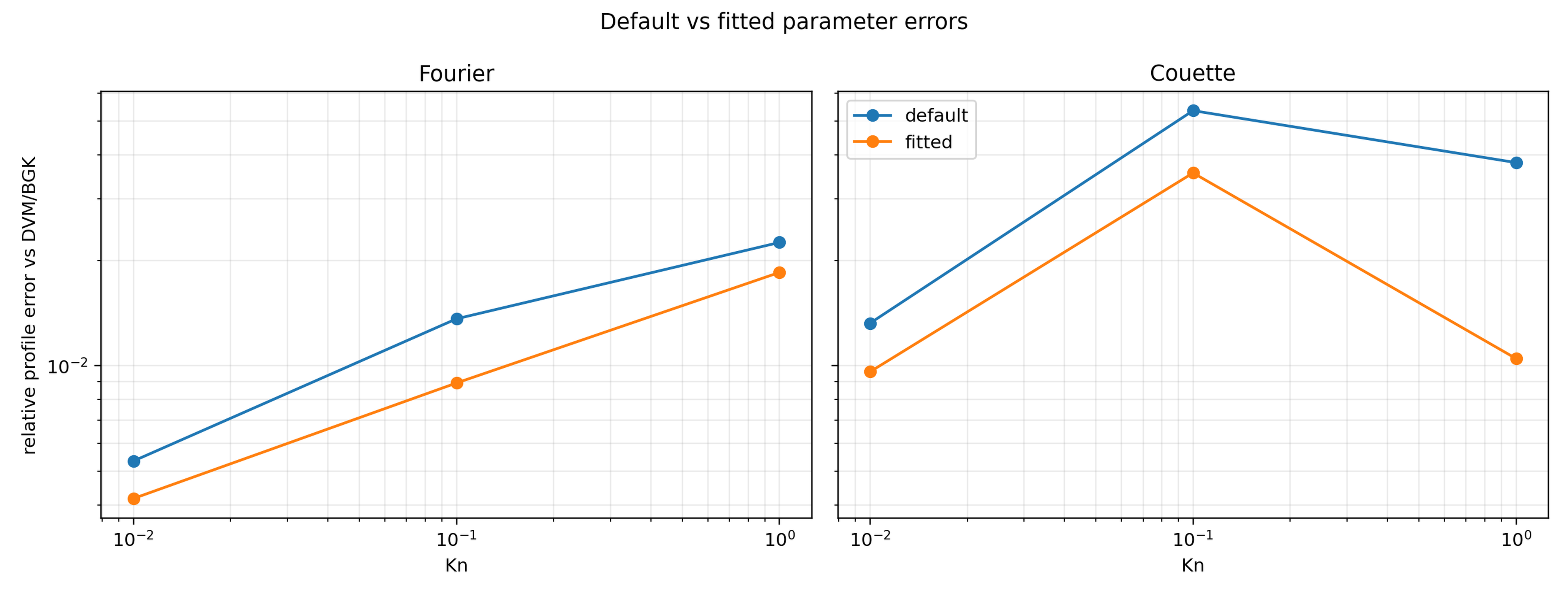}
\caption{Relative profile errors for default and calibrated log-Gaussian parameters.}
\label{fig:default_fitted_errors}
\end{figure}

The calibrated pair should be interpreted as a DVM/BGK-informed recommendation for the present planar wall-bounded benchmark family, not as a universal optimum.


\subsection{Preliminary shock-like indicator diagnostic}
This subsection provides a preliminary analytic sensor check; the reduced DVM/BGK
shock-layer calculation below provides the stronger kinetic activation diagnostic. The
test is not intended as a shock-structure validation against DSMC or DVM data. Instead,
it checks whether the proposed local indicator activates the kinetic branch inside a
shock-like high-gradient layer even when the global Knudsen number is small.

A smooth ideal-gas normal-shock-like profile is constructed from Rankine--Hugoniot
upstream and downstream states using a hyperbolic-tangent transition layer with
thickness proportional to the global mean free path. The local indicators
\[
K_\rho=\lambda\frac{|\partial_x\rho|}{\rho},\qquad
K_T=\lambda\frac{|\partial_xT|}{T},\qquad
K_u=\lambda\frac{|\partial_xu|}{|u|+a},
\]
are evaluated across the layer and combined as \(K=\max(K_\rho,K_T,K_u)\). For
\(M_1=3\) and global \(Kn=10^{-3}\), the maximum local indicator is
\(K_{\max}=5.43\times10^{-2}\), giving maximum ballistic activation
\(W_f=0.271\) with the default parameters and \(W_f=0.594\) with the calibrated
parameters. For \(M_1=5\) and global \(Kn=10^{-3}\), the corresponding values are
\(K_{\max}=7.64\times10^{-2}\), \(W_f=0.394\), and \(W_f=0.646\). Thus even this
analytic preliminary test shows that a gradient-based local indicator can activate the
kinetic branch inside shock-like layers when the global rarefaction parameter is small.

\subsection{Leave-one-case-out parameter robustness check}

The calibrated pair \(K_0=0.03,\sigma=2.5\) was obtained by minimizing the mean DVM/BGK profile error over all six Fourier and Couette benchmark cases. To test whether this choice is merely an artifact of fitting all cases simultaneously, we perform a leave-one-case-out internal robustness check across this closely related planar benchmark family. Each physical case, defined by the flow type and Knudsen number, is held out in turn. The parameters are selected by minimizing the mean profile error over the remaining five cases and are then evaluated on the held-out case.

Table~\ref{tab:holdout_cv} summarizes the result. The training-selected parameters reduce the mean held-out error from \(2.44\times10^{-2}\) for the default pair to \(2.10\times10^{-2}\). The median selected parameters are again \(K_0=0.03,\sigma=2.5\), matching the full-data optimum. This indicates internal stability of the selected transition scale within the present benchmark family.

\begin{table}[t]
\centering
\caption{Leave-one-case-out parameter robustness check. Parameters are selected by minimizing the mean profile error over the remaining five DVM/BGK benchmark cases and then evaluated on the held-out case.}
\label{tab:holdout_cv}
\begin{tabular}{lcccc}
\toprule
Choice & Mean held-out error & Median selected $K_0$ & Median selected $\sigma$ & Cases \\
\midrule
default & 2.4368e-02 & 0.1 & 1.0 & 6 \\
training-selected & 2.1048e-02 & 0.03 & 2.5 & 6 \\
full-data best & 1.4533e-02 & 0.03 & 2.5 & 6 \\
\bottomrule
\end{tabular}
\end{table}

\begin{figure}[t]
\centering
\includegraphics[width=0.88\textwidth]{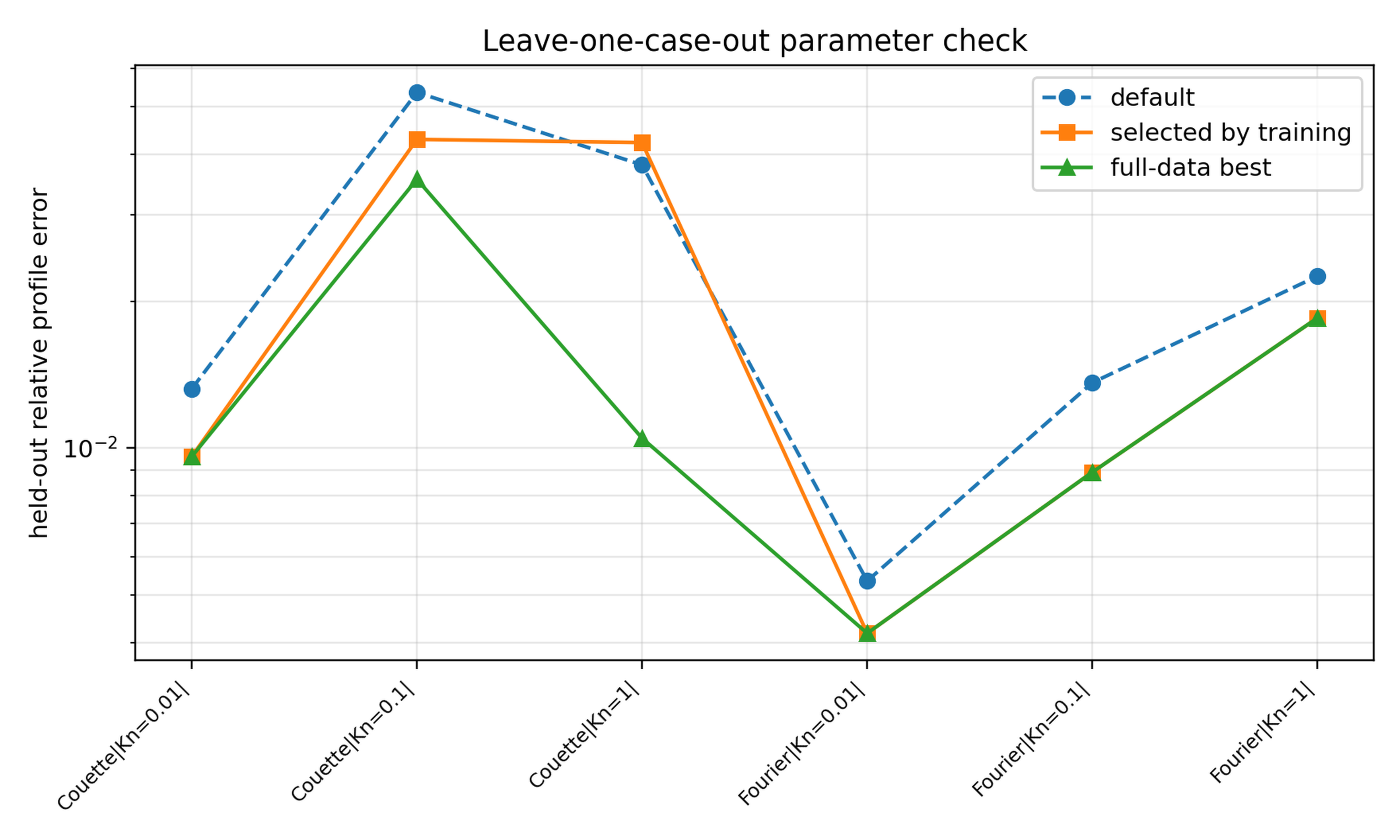}
\caption{Leave-one-case-out internal robustness check. Parameters selected on the remaining five closely related planar cases improve the mean held-out error relative to the default pair, but the Couette \(Kn=1\) case shows that the small calibration set does not establish independent validation or universal parameter transfer.}
\label{fig:holdout_cv}
\end{figure}

Figure~\ref{fig:holdout_cv} shows the held-out errors case by case. The check is intentionally conservative and internal: it probes sensitivity within six closely related planar Fourier/Couette cases, but it does not replace independent validation across different flow classes, external datasets, DSMC, UGKS, or experimental data. It also exposes a limitation. In the held-out Couette \(Kn=1\) case, the parameters selected without that case are less accurate than the default, while the full-data calibrated pair remains accurate. This indicates that the available six-case calibration set is still small and that broader external validation is required before treating the parameters as universal.

\FloatBarrier

\label{sec:additional_diagnostic_suite}
For the remaining component and shock-layer diagnostics, the preceding sections established the log-Gaussian limiter as a proof-of-concept continuum--ballistic
blending mechanism and calibrated the transition parameters against planar DVM/BGK
Fourier and Couette benchmarks. We next examine five diagnostic components. These tests separate the roles of the BGK relaxation coefficient, non-equilibrium moment accuracy, parameter robustness, one-dimensional kinetic shock-layer activation, and two-dimensional curved shock-layer activation. They are diagnostic tests only and do not constitute external DSMC or full hypersonic blunt-body validation, but they substantially strengthen the evidence that the local rarefaction indicator activates in physically relevant high-gradient regions.

\subsection{Verification of the BGK time-averaging coefficient}

The BGK-corrected kinetic flux uses
\[
F_n^K = \alpha F_n^{FM} + (1-\alpha)F_n^{EQ},
\qquad
\alpha = \frac{\tau}{\Delta t}\left(1-e^{-\Delta t/\tau}\right).
\]
This coefficient is not an empirical damping factor. It is the exact time-averaging coefficient for exponential relaxation from an upwind free-transport distribution toward an interface equilibrium distribution over one time step.

To verify this interpretation, we compare the model flux with a direct discrete-velocity time-averaged BGK relaxation calculation for several left--right Maxwellian interface states. The maximum relative flux error is approximately machine precision, about \(4\times 10^{-16}\), across the tested cases. This confirms the internal consistency of the relaxation coefficient used in the kinetic branch. The result should not be interpreted as a full UGKS asymptotic-preserving proof; it verifies the local interface relaxation coefficient only.

\begin{figure}[t]
\centering
\includegraphics[width=0.48\textwidth]{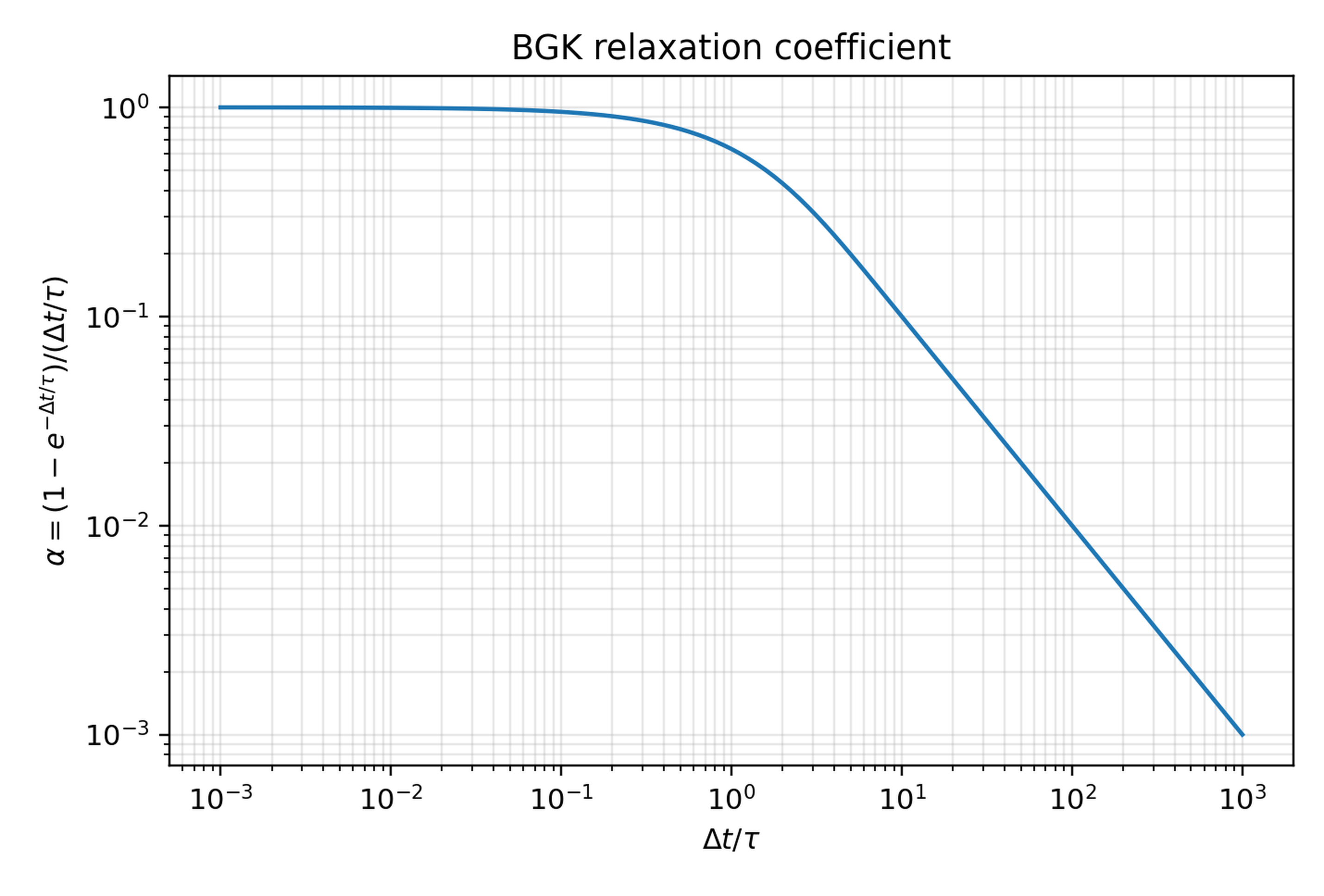}
\includegraphics[width=0.48\textwidth]{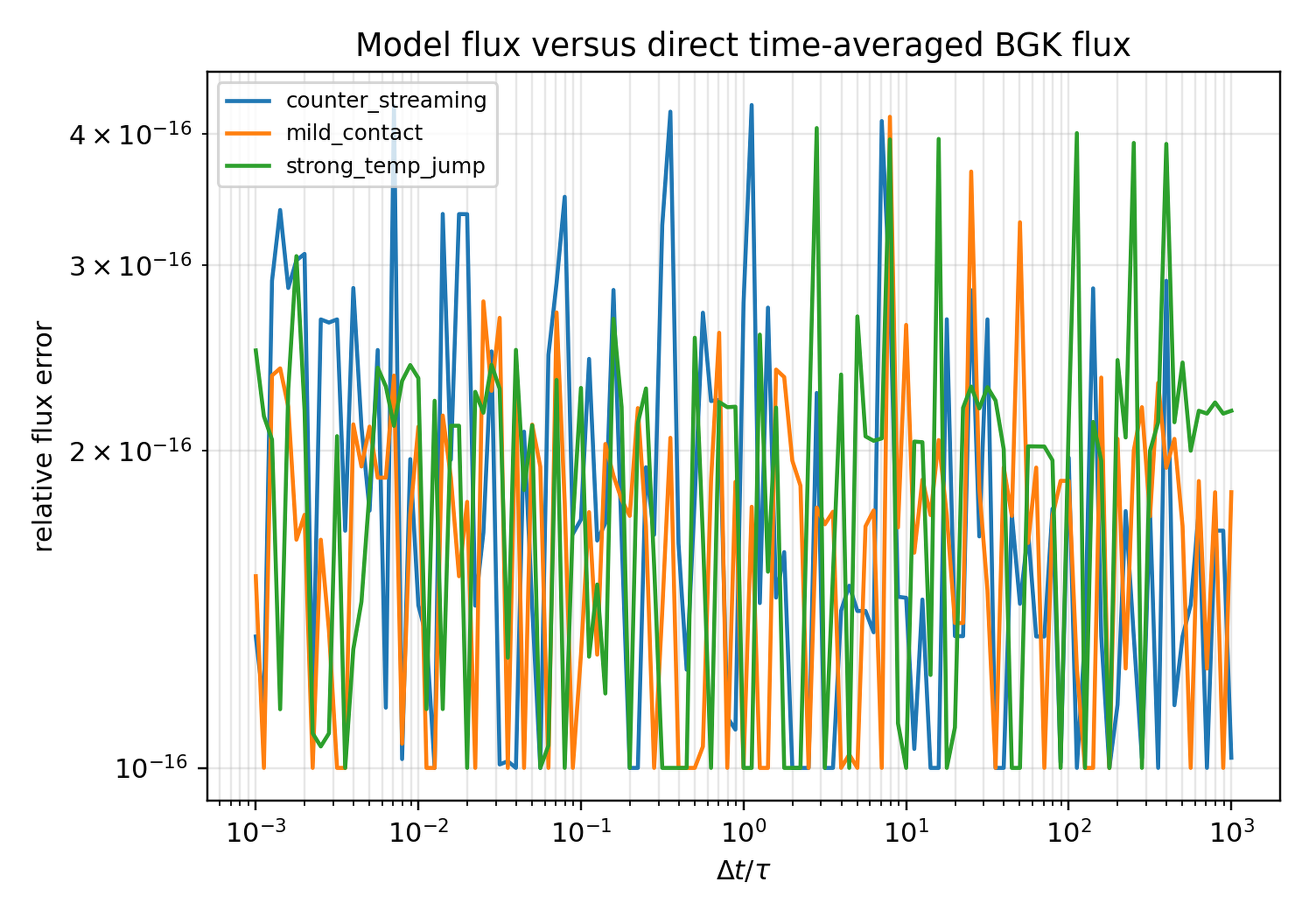}
\caption{BGK relaxation coefficient diagnostic. Left: the relaxation weight \(\alpha\) as a function of \(\Delta t/\tau\). Right: relative error between the blended model flux and the direct time-averaged BGK relaxation flux.}
\label{fig:bgk_relaxation_validation}
\end{figure}

\subsection{Non-equilibrium heat-flux and shear-stress diagnostics}

The profile-based assessment above does not directly examine non-equilibrium fluxes and moments. We therefore compare heat-flux and shear-stress discrepancies against the DVM/BGK reference data within a reference-informed diagnostic construction.

For the Fourier heat-flux diagnostic, the predictor is constructed as
\begin{equation}
q_{model}=(1-W_f)q_{NSF}+W_f q_{jump},
\end{equation}
with
\begin{equation}
q_{NSF}=-\kappa_{eff}\frac{\Delta T}{H},
\qquad
q_{jump}=-\kappa_{eff}\frac{\Delta T}{H+2B_TKH},
\qquad
\kappa_{eff}=-\frac{q_{DVM}}{\langle dT_{DVM}/dy\rangle}.
\end{equation}
For the Couette shear-stress diagnostic, the model stress is computed from the reduced velocity profile as
\begin{equation}
\tau_{model}=\mu_{eff}\left\langle\frac{du_{model}}{dy}\right\rangle,
\qquad
\mu_{eff}=\frac{\tau_{DVM}}{\langle du_{DVM}/dy\rangle}.
\end{equation}
Here \(\langle\cdot\rangle\) denotes the interior average used in the post-processing scripts, excluding ten grid points near each wall. These effective coefficients are extracted from the DVM/BGK references for diagnostic comparison; they should not be interpreted as independent first-principles transport closures.

\begin{table}[t]
\centering
\caption{DVM/BGK profile, flux, and moment diagnostics. Errors are relative to DVM/BGK reference values. Within the reference-informed diagnostics, rarefaction branches reduce the reported heat-flux and shear-stress discrepancies relative to NSF, while the log-Gaussian limiter is not always the best moment closure.}
\label{tab:flux_moment_validation}
\begin{tabular}{lllccc}
\toprule
Benchmark & Moment & Model & Mean error & Max error & Cases \\
\midrule
Fourier & heat\_flux & nsf & 8.0438e-01 & 2.0450e+00 & 3 \\
Fourier & heat\_flux & hard & 1.8286e-01 & 4.3190e-01 & 3 \\
Fourier & heat\_flux & algebraic & 1.9120e-01 & 4.0738e-01 & 3 \\
Fourier & heat\_flux & log & 1.9057e-01 & 4.0552e-01 & 3 \\
Fourier & heat\_flux & jump & 1.7007e-01 & 4.3190e-01 & 3 \\
Couette & velocity\_profile & nsf & 1.5648e-01 & 3.5363e-01 & 3 \\
Couette & velocity\_profile & hard & 2.1035e-02 & 4.2230e-02 & 3 \\
Couette & velocity\_profile & algebraic & 3.5018e-02 & 5.3505e-02 & 3 \\
Couette & velocity\_profile & log & 3.4916e-02 & 5.3505e-02 & 3 \\
Couette & velocity\_profile & slip & 1.7297e-02 & 4.2230e-02 & 3 \\
Couette & shear\_stress & nsf & 6.3536e-01 & 1.6396e+00 & 3 \\
Couette & shear\_stress & hard & 7.7912e-02 & 1.9819e-01 & 3 \\
Couette & shear\_stress & algebraic & 1.1027e-01 & 1.8000e-01 & 3 \\
Couette & shear\_stress & log & 1.0980e-01 & 1.7862e-01 & 3 \\
Couette & shear\_stress & slip & 7.0236e-02 & 1.9819e-01 & 3 \\
\bottomrule
\end{tabular}
\end{table}

Within this reference-informed diagnostic construction, the rarefaction branches reduce the reported heat-flux and shear-stress discrepancies relative to NSF. For Fourier heat flux, the NSF mean relative error is about \(8.04\times 10^{-1}\), whereas the log-Gaussian limiter gives about \(1.91\times 10^{-1}\). For Couette shear stress, the NSF mean relative error is about \(6.35\times 10^{-1}\), whereas the log-Gaussian limiter gives about \(1.10\times 10^{-1}\).

However, the log-Gaussian limiter is not always the best moment closure. Dedicated jump or slip branches can give smaller moment errors. This distinction is important: the log-Gaussian weighting changes the transition between continuum and rarefied branches, but the reported heat-flux and shear-stress discrepancies remain tied to the reference-informed diagnostic construction and to the rarefied-branch closure used. Because these moment comparisons are constructed from the existing DVM/BGK benchmark metrics and use effective transport scaling in the rarefied branches, they should be interpreted as moment-error diagnostics rather than independent first-principles closure validation.

\subsection{Extended parameter robustness}

The calibrated parameter pair \(K_0=0.03,\sigma=2.5\) was obtained from a small DVM/BGK benchmark family. To test whether this result is overly dependent on a single aggregation procedure, we perform a case-weighted ranking, bootstrap resampling over physical cases, and per-case optimum analysis.

\begin{table}[t]
\centering
\caption{Extended parameter calibration robustness analysis based on the available DVM/BGK benchmark family. Case-weighted errors average each physical case equally. Bootstrap fractions are obtained by resampling cases with replacement.}
\label{tab:exp18_parameter_robustness}
\begin{tabular}{lcc}
\toprule
Quantity & Value & Comment \\
\midrule
Case-weighted best $K_0$ & 0.03 & minimum mean error \\
Case-weighted best $\sigma$ & 2.5 & minimum mean error \\
Best mean error & 1.4533e-02 & case-weighted \\
Default mean error & 2.4368e-02 & $K_0=0.1,\sigma=1.0$ \\
Fitted mean error & 1.4533e-02 & $K_0=0.03,\sigma=2.5$ \\
Top bootstrap pair & $K_0=0.05,\sigma=2.5$ & selected most often \\
Top bootstrap fraction & 0.374 & resampled cases \\
\bottomrule
\end{tabular}
\end{table}

The case-weighted optimum remains \(K_0=0.03,\sigma=2.5\), with mean profile error \(1.4533\times10^{-2}\). Nearby broad-transition choices, especially with \(\sigma=2.5\), have similar errors. This indicates a broad, relatively flat optimum valley rather than a sharply isolated parameter point.

Bootstrap resampling gives a more cautious picture. The most frequently selected pairs are \(K_0=0.05,\sigma=2.5\), \(K_0=0.01,\sigma=0.25\), and \(K_0=0.03,\sigma=2.5\). Thus the calibrated pair is best interpreted as a representative global compromise for the present benchmark family, not as a universal constant. Per-case optima also vary with flow type and Knudsen number, suggesting that state-dependent or regime-adaptive transition parameters may be useful in future work.

\begin{figure}[t]
\centering
\includegraphics[width=0.48\textwidth]{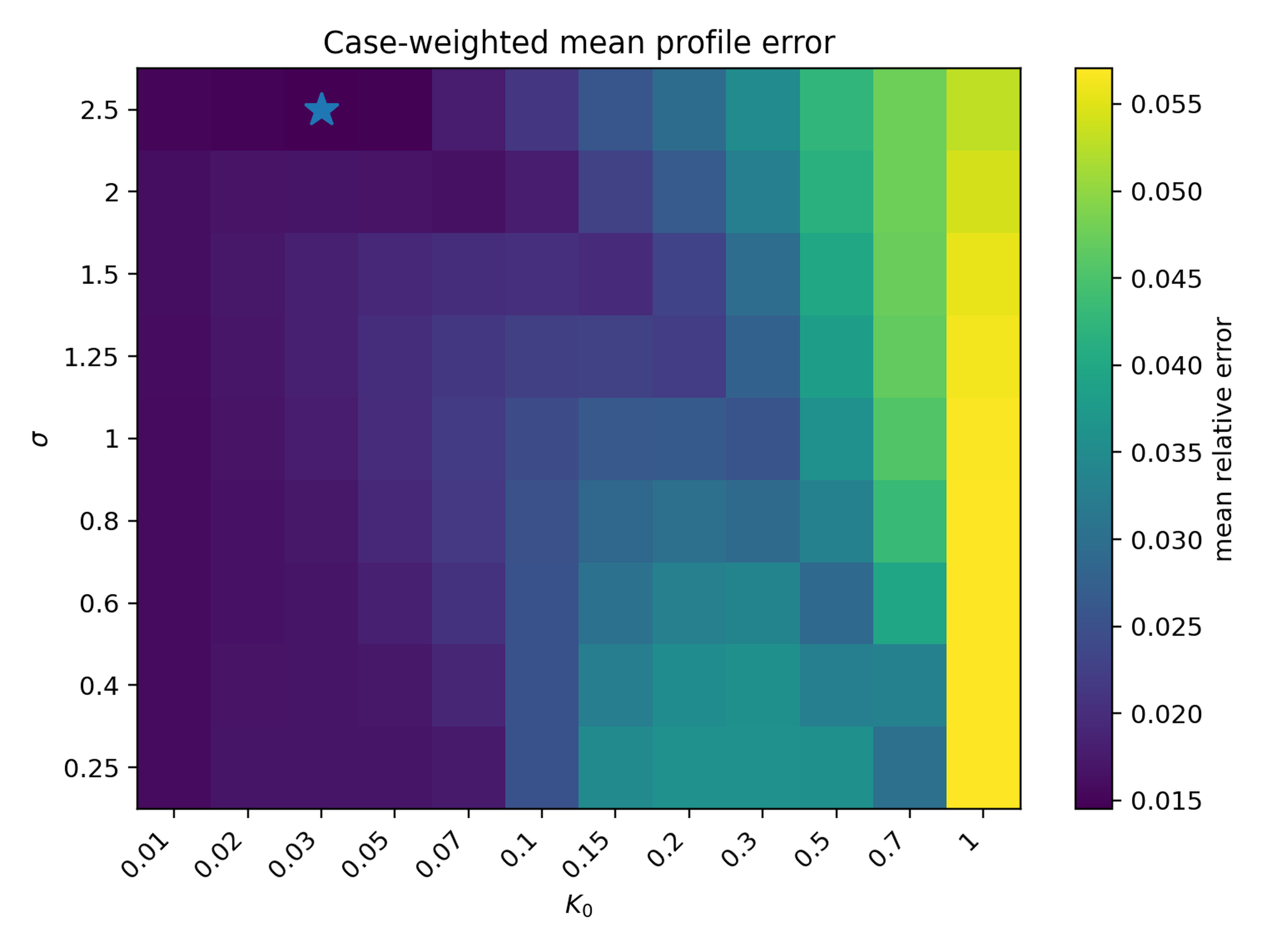}
\includegraphics[width=0.48\textwidth]{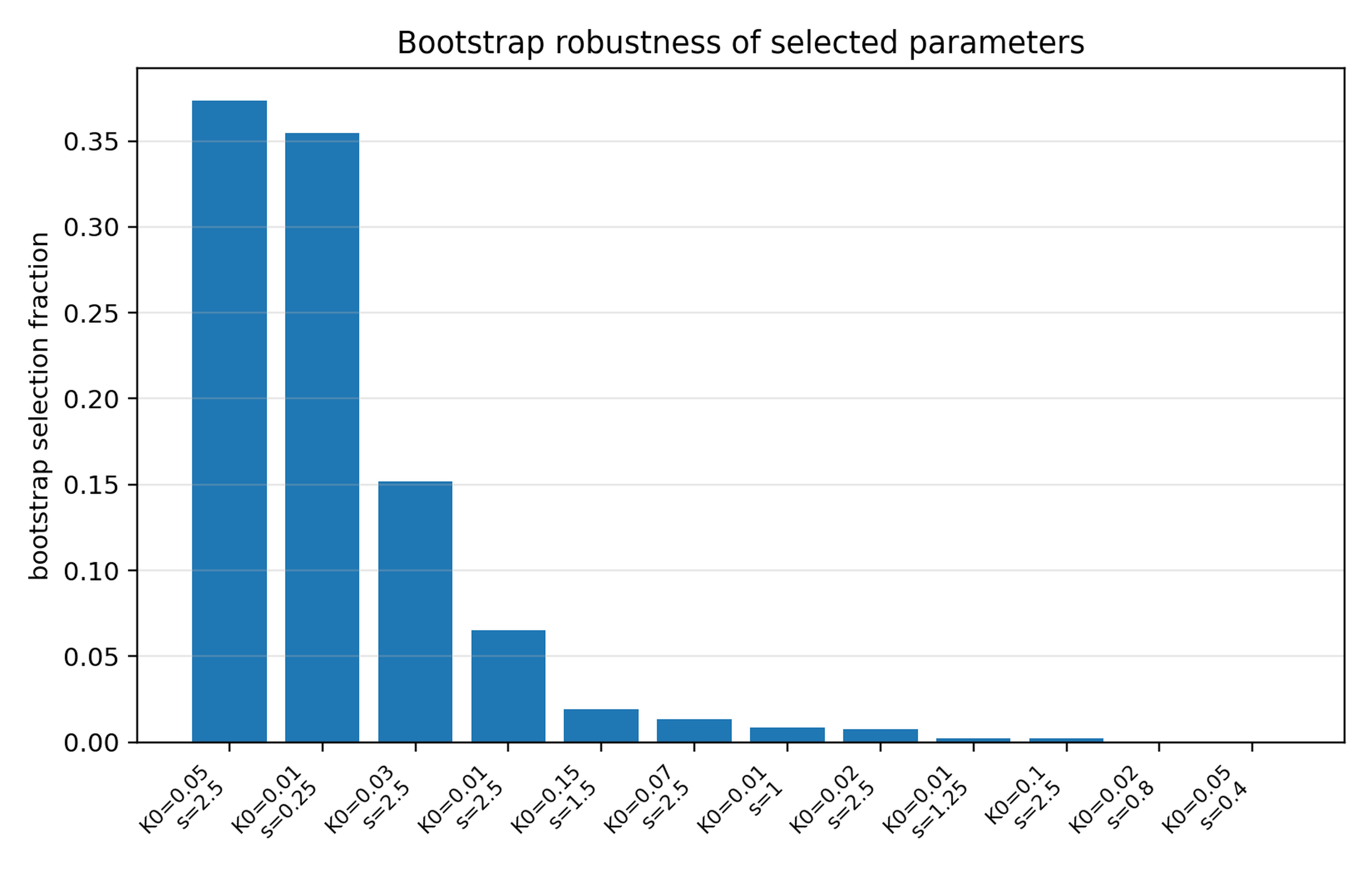}
\caption{Extended parameter robustness analysis. Left: case-weighted mean profile error over the scanned \((K_0,\sigma)\) grid. Right: bootstrap selection fractions obtained by resampling the available DVM/BGK cases.}
\label{fig:parameter_robustness}
\end{figure}

\begin{figure}[t]
\centering
\includegraphics[width=0.82\textwidth]{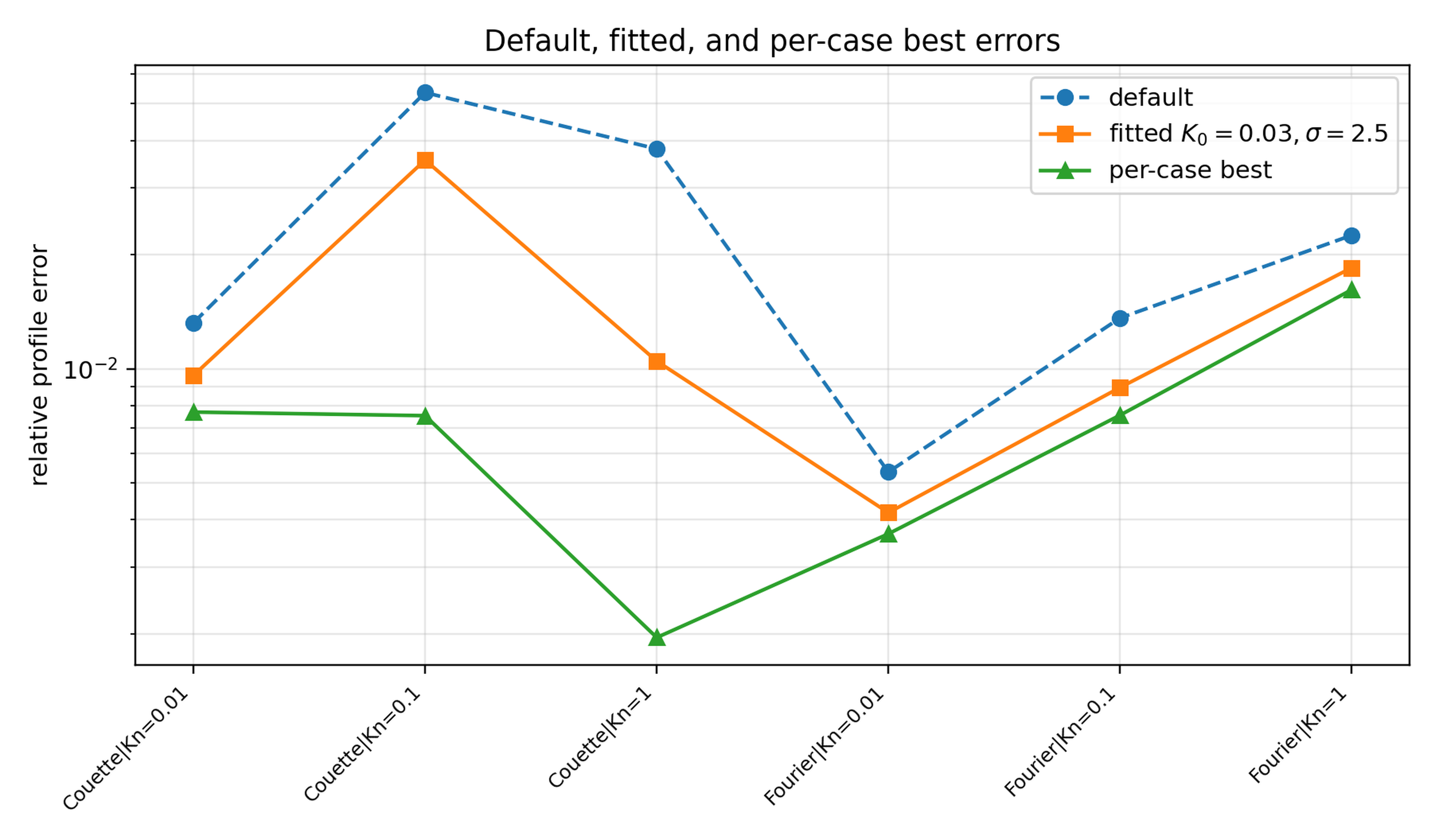}
\caption{Default, fitted, and per-case best profile errors. The fitted pair is a useful global compromise, while individual cases can prefer different transition parameters.}
\label{fig:per_case_parameters}
\end{figure}

\subsection{Reduced one-velocity DVM/BGK shock-layer activation}
\label{sec:reduced_dvm_bgk_shock}

We next test the local rarefaction indicator in a kinetic shock layer. The diagnostic uses a reduced one-dimensional DVM/BGK equation,
\[
\partial_t f + v\,\partial_x f = \frac{M[f]-f}{\tau},
\]
with upstream and downstream Maxwellian states based on Rankine--Hugoniot normal-shock relations. The model uses one velocity dimension and is intended as a kinetic shock-layer diagnostic. It is not a DSMC validation of monatomic-gas shock structure and not a full hypersonic blunt-body solver.

\begin{table}[t]
\centering
\caption{Reduced one-velocity DVM/BGK shock-layer diagnostic. The test uses a one-dimensional BGK kinetic model to examine activation of the local rarefaction indicator inside shock layers. It is a kinetic diagnostic, not a DSMC validation of monatomic-gas shock structure.}
\label{tab:exp19_normal_shock}
\begin{tabular}{cccccccc}
\toprule
$M_1$ & Kn & $\max K_{\rm local}$ & $\max W_f$ default & $\max W_f$ fitted & Thickness & Residual & Conv. \\
\midrule
2.0 & 0.03 & 0.1971 & 0.7513 & 0.7743 & 0.2372 & 8.99e-07 & yes \\
3.0 & 0.03 & 0.2445 & 0.8144 & 0.7993 & 0.2776 & 9.47e-07 & yes \\
3.0 & 0.10 & 0.2881 & 0.8550 & 0.8172 & 0.8298 & 9.96e-07 & yes \\
5.0 & 0.03 & 0.1995 & 0.7551 & 0.7757 & 0.4105 & 1.00e-06 & yes \\
\bottomrule
\end{tabular}
\end{table}

All cases reach the \(10^{-6}\) relative-update tolerance. Across Mach 2--5 cases, the maximum local rarefaction indicator reaches approximately \(K_{\rm local}=0.20\)--\(0.29\), producing maximum kinetic weights of about \(W_f=0.75\)--\(0.86\). At \(M_1=3\), increasing Kn from 0.03 to 0.10 increases the measured shock thickness from about 0.28 to about 0.83. These results support the use of a local gradient-based rarefaction indicator in shock-layer regions.

\begin{figure}[t]
\centering
\includegraphics[width=0.82\textwidth]{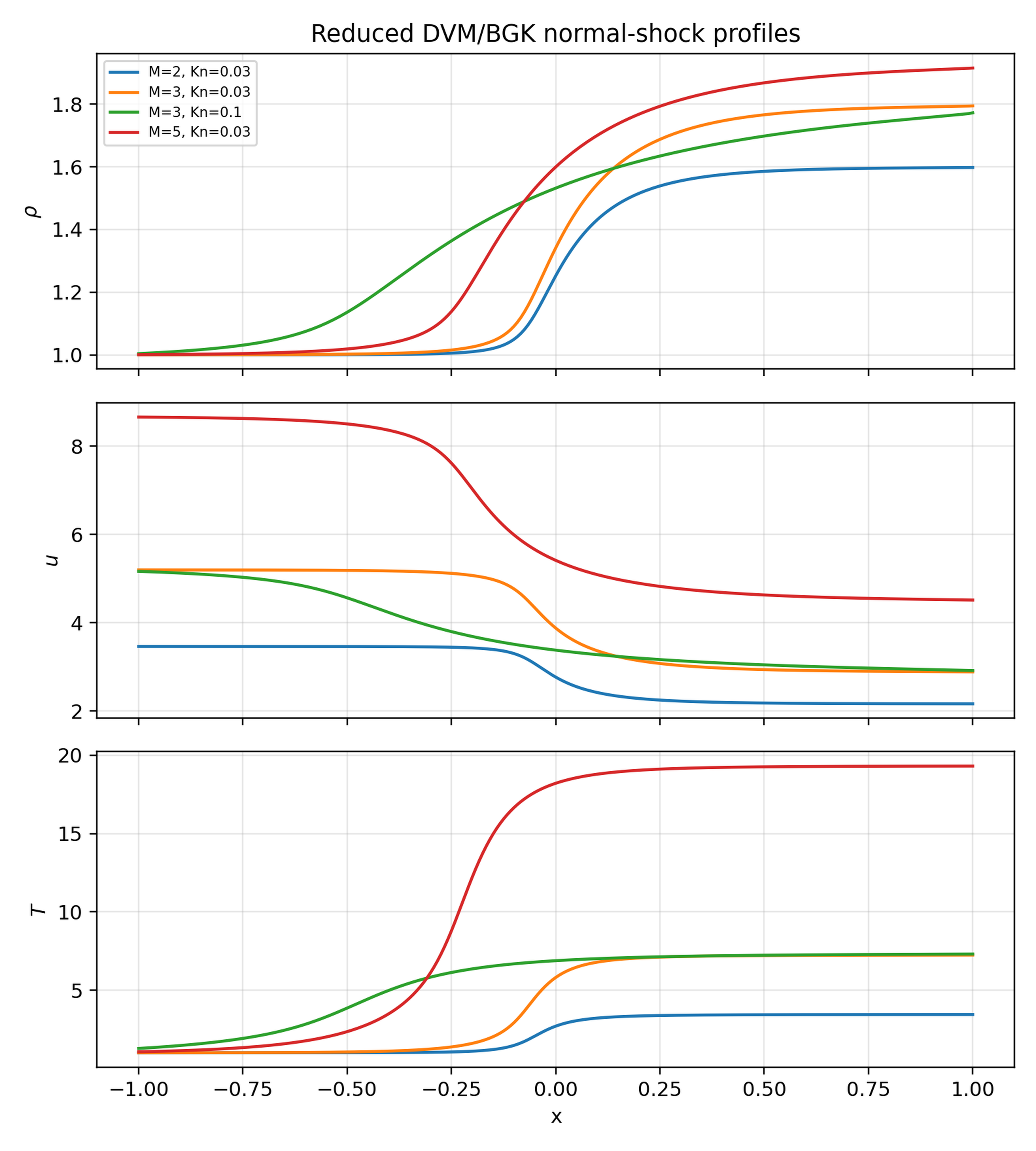}
\caption{Reduced one-velocity DVM/BGK shock-layer profiles for density, velocity, and temperature.}
\label{fig:normal_shock_profiles}
\end{figure}

\begin{figure}[t]
\centering
\includegraphics[width=0.82\textwidth]{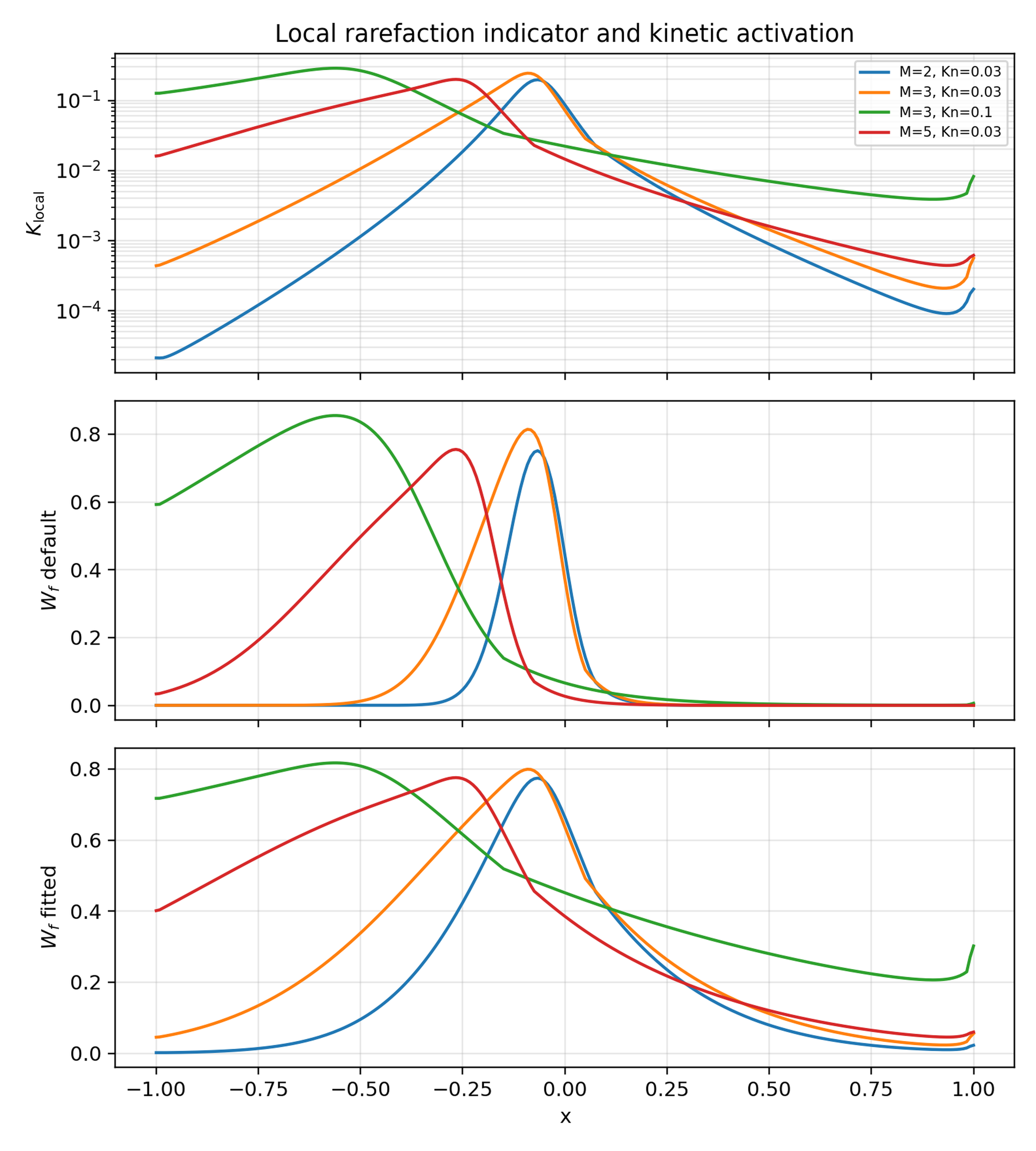}
\caption{Local rarefaction indicator and kinetic activation weights in the reduced one-velocity DVM/BGK shock-layer diagnostic.}
\label{fig:normal_shock_indicator}
\end{figure}

\subsection{Two-dimensional curved shock-layer activation diagnostic}

Finally, we evaluate the indicator in a two-dimensional curved bow-shock geometry. The field is a smooth analytic blunt-body shock-layer construction using normal-shock jump relations and a curved bow-shock profile. This test is not a full 2D hypersonic CFD calculation and is not a DSMC validation. Its purpose is to test whether the local rarefaction indicator activates in a curved shock-layer and stagnation-region geometry.

\begin{table}[t]
\centering
\caption{Two-dimensional blunt-body shock-layer diagnostic. The field is a smooth analytic bow-shock construction used to test the local rarefaction indicator in curved shock geometry. It is not a full hypersonic CFD or DSMC validation.}
\label{tab:exp20_blunt_body}
\begin{tabular}{ccccccc}
\toprule
$M_\infty$ & Kn & $\max K_{\rm local}$ & $\max W_f$ def. & $\max W_f$ fit. & Area $K>0.1$ & Area $W_f>0.5$ fit. \\
\midrule
5.0 & 0.010 & 0.0943 & 0.4764 & 0.6765 & 0.0000 & 0.0418 \\
5.0 & 0.030 & 0.2818 & 0.8499 & 0.8149 & 0.0053 & 0.1152 \\
5.0 & 0.100 & 0.9569 & 0.9880 & 0.9170 & 0.0053 & 0.3918 \\
8.0 & 0.030 & 0.2594 & 0.8297 & 0.8059 & 0.0261 & 0.1488 \\
\bottomrule
\end{tabular}
\end{table}

For \(M_\infty=5\), increasing Kn from 0.01 to 0.10 raises the maximum local indicator from about 0.094 to about 0.957. The fitted kinetic-activation area, measured by the fraction of the fluid domain with \(W_f>0.5\), increases from about 4.2\% to about 39.2\%. At \(M_\infty=8\), \(Kn=0.03\), the high-gradient region with \(K_{\rm local}>0.1\) is larger than in the corresponding \(M_\infty=5\) case. These results show that the gradient-based local indicator is sensitive to curved shock-layer geometry.

\begin{figure}[t]
\centering
\includegraphics[width=0.9\textwidth]{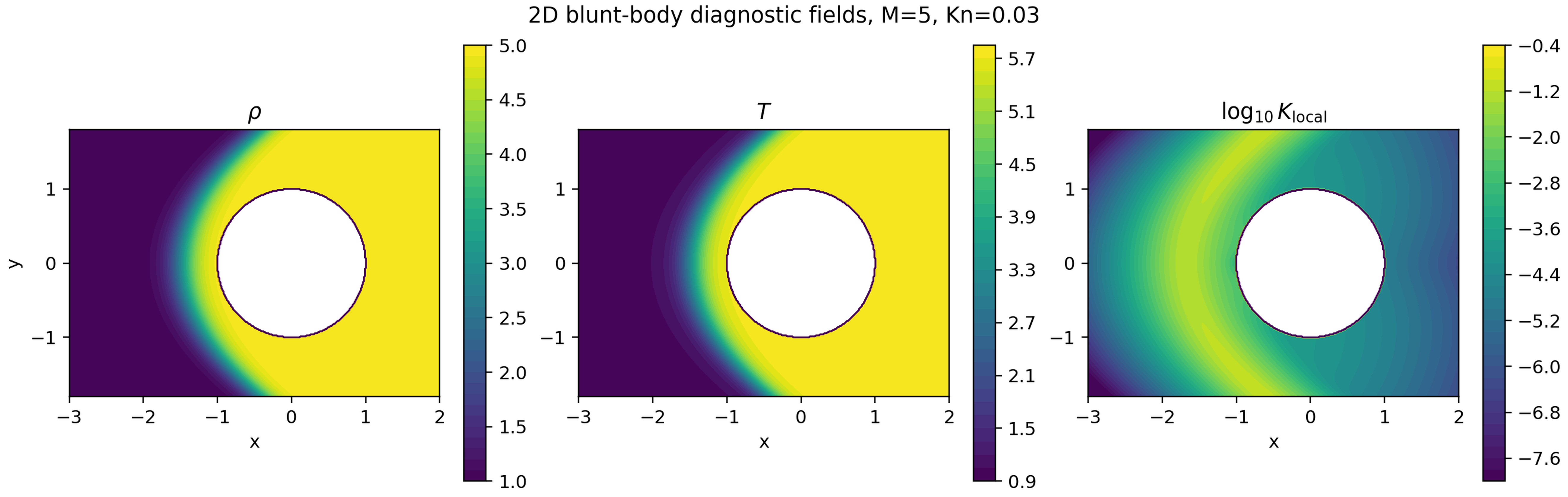}
\caption{Two-dimensional blunt-body diagnostic fields for density, temperature, and local rarefaction indicator.}
\label{fig:blunt_body_fields}
\end{figure}

\begin{figure}[t]
\centering
\includegraphics[width=0.82\textwidth]{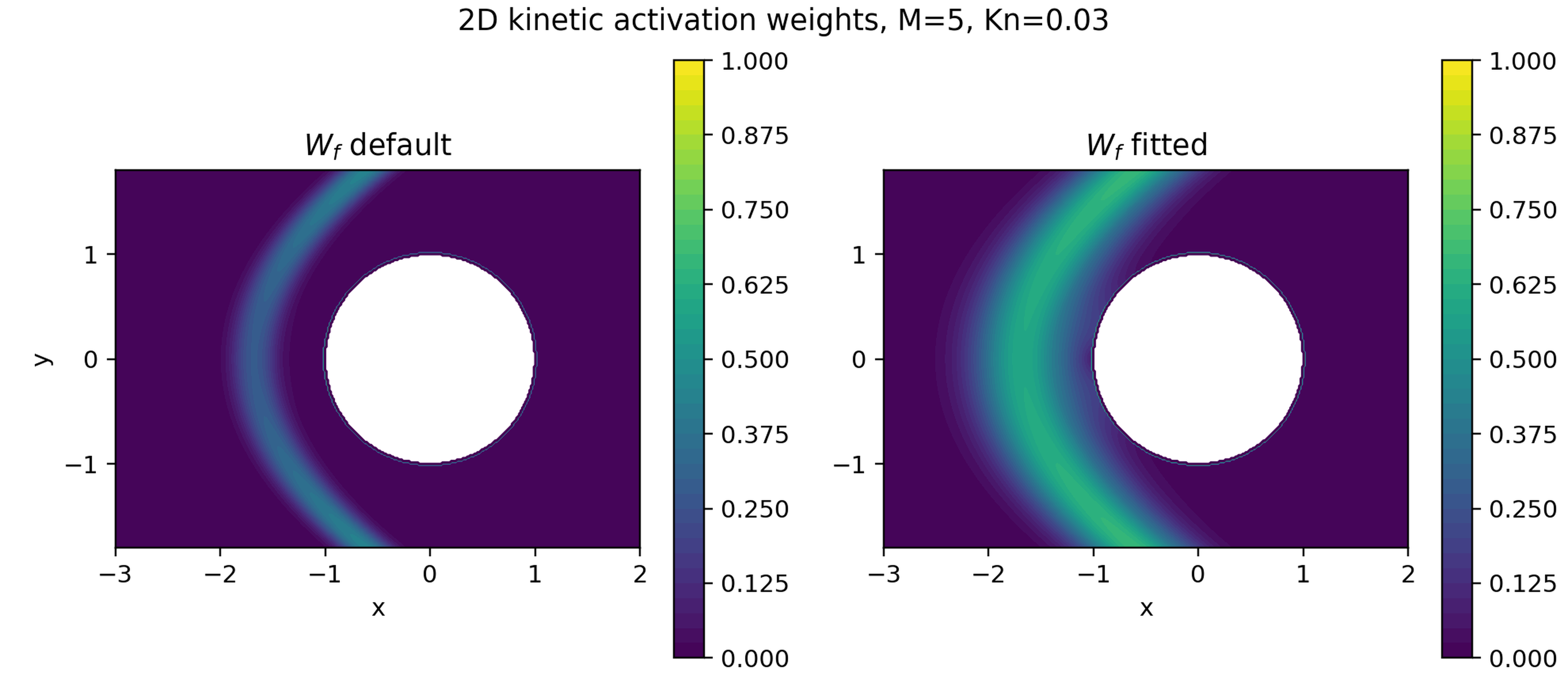}
\caption{Default and fitted log-Gaussian kinetic activation weights in the two-dimensional blunt-body shock-layer diagnostic.}
\label{fig:blunt_body_weights}
\end{figure}

\subsection{Summary of diagnostics}

The additional diagnostics clarify the current status of the method. The BGK coefficient has an exact local relaxation interpretation. Within the reference-informed diagnostic construction, rarefaction branches reduce the reported heat-flux and shear-stress discrepancies relative to NSF, but moment accuracy remains limited by the rarefied-branch closure. The calibrated transition parameters remain optimal under case-weighted aggregation, but bootstrap and per-case analyses show non-negligible data-set dependence. The local rarefaction indicator activates strongly in both reduced kinetic normal-shock layers and two-dimensional curved shock-layer geometries.

Together, these results extend the assessment beyond planar calibration to non-equilibrium moments, parameter sensitivity, and shock-layer activation. At the same time, they do not remove the need for future DSMC, UGKS, or experimental validation of full hypersonic rarefied flows.

\FloatBarrier

\section{Discussion}

The proposed limiter differs from a purely empirical switch because the weights are defined in
logarithmic Knudsen space and can be interpreted as a probability partition between collisional
and ballistic transport. This interpretation is useful in rarefied and hypersonic gas dynamics,
where the relevant transition is not controlled by a single global Knudsen number but by local
mean-free-path-to-gradient-length ratios. The use of a face-local indicator allows the kinetic
branch to activate in shocks, wall layers, and other high-gradient regions even when the global
rarefaction level is small.

The planar DVM/BGK Fourier and Couette data are used both to construct reference profiles
for the reduced comparisons and to calibrate \(K_0\) and \(\sigma\). Within this in-sample
calibration setting, rarefaction corrections improve macroscopic profiles relative to NSF. The
calibrated pair \(K_0=0.03,\sigma=2.5\) performs better than the default \(K_0=0.1,\sigma=1.0\)
for the present benchmark set, suggesting that wall-induced kinetic effects can become visible
before the conventional \(Kn=0.1\) threshold and that a broader transition in logarithmic
Knudsen space is beneficial for these flows. The extended case-weighted and bootstrap analyses
refine this conclusion: the calibrated pair remains the best case-weighted compromise, but
bootstrap resampling and per-case optima reveal a broad parameter valley and non-negligible
data-set dependence. The parameters should therefore be interpreted as DVM/BGK-informed
calibration values for the current benchmark family, not as universal constants or independently
validated coefficients.

The additional diagnostics clarify which parts of the framework are robust and which parts
remain closure-dependent. The BGK relaxation coefficient has an exact local time-averaging
interpretation for exponential relaxation from a free-transport interface state toward an
equilibrium state. Within the reference-informed heat-flux and shear-stress diagnostics, rarefaction branches reduce the reported discrepancies relative to NSF, but the best moment accuracy is often obtained by dedicated slip or jump branches. Thus the log-Gaussian weighting improves
the transition mechanism, while accurate heat fluxes and stresses still require sufficiently
accurate rarefied-branch closures.

The shock-layer diagnostics support the use of a local rarefaction indicator. In both the
smooth normal-shock-like diagnostic and the reduced one-velocity DVM/BGK shock-layer calculation,
the local indicator reaches values large enough to activate the kinetic branch inside the shock
layer. The two-dimensional blunt-body diagnostic further shows that the same mechanism
activates in curved bow-shock geometry and stagnation-region-like fields. These tests strengthen
the interpretation of the method as a local continuum--kinetic activation framework.

The present study still has important limitations. The reduced one-velocity DVM/BGK shock-layer calculation is
not a DSMC validation of monatomic shock structure, and the two-dimensional blunt-body case
is an analytic diagnostic field rather than a full hypersonic CFD computation. Strong shocks,
highly non-Maxwellian distributions, separated multidimensional flows, and realistic gas effects
may require different indicators, adaptive transition parameters, or more accurate kinetic
branches. The present framework is modular and can be combined with more sophisticated
continuum fluxes, kinetic fluxes, DVM/DSMC-informed closures, or moment-based rarefied
branches.

\section{Conclusions}

A log-Gaussian scale-space limiter has been proposed for hybrid continuum--ballistic gas
dynamics. By treating the local Knudsen number as a stochastic scale ratio, the continuum
and ballistic weights are defined as complementary cumulative probabilities in logarithmic
Knudsen space. The resulting weights recover the continuum and free-molecular limits and
super-algebraically suppress asymptotically invalid correction branches.

The limiter is incorporated into a conservative interface flux blending continuum and kinetic
contributions. Reduced one-dimensional closure/profile comparisons against DVM/BGK Fourier and Couette data show that log-Gaussian weighting of NSF and jump/slip-corrected branches improves the tested macroscopic profiles relative to NSF. Because the same six profiles are used to calibrate \(K_0\) and \(\sigma\), the resulting approximately 40\% reduction in combined mean profile error is an in-sample calibration result for the reduced profile model. It is not an independent validation result and does not numerically validate the proposed finite-volume face flux in a complete solver.

Additional diagnostics strengthen the evidence base of the study. A direct interface
relaxation test verifies the BGK coefficient as the exact time-averaging factor for exponential
BGK relaxation. Within the reference-informed heat-flux and shear-stress diagnostics, rarefaction branches reduce the reported discrepancies relative to NSF, while also demonstrating that moment accuracy remains limited by the rarefied-branch closure. Case-weighted and bootstrap parameter analyses show
that the calibrated pair is a useful global compromise for the present DVM/BGK benchmark
family, but not a universal optimum.

Shock-layer diagnostics show that the proposed local rarefaction indicator activates the kinetic
branch in high-gradient regions. This activation is observed in a smooth shock-like diagnostic,
a reduced one-velocity DVM/BGK shock-layer calculation, and a two-dimensional blunt-body bow-shock
diagnostic. These results support the use of a local gradient-based indicator rather than a
purely global Knudsen-number switch.

The present study is a proof-of-concept framework, in-sample calibration study, and diagnostic assessment rather than a replacement for kinetic solvers or a completed finite-volume solver validation. The hybrid face flux is formulated for future implementation in finite-volume and gas-kinetic solvers, but it has not been advanced in a complete production solver in the present study. Future work will include DSMC, UGKS, and iterative kinetic-solver comparisons \cite{su_zhu_wang_zhang_wu_gsis}, full normal-shock structure validation, multidimensional hypersonic flows,
improved non-Maxwellian kinetic branches, adaptive transition parameters, and analysis of
positivity, entropy stability, and asymptotic-preserving behavior.


\end{document}